\documentclass[a4paper]{article}

\usepackage{geometry}

\usepackage{palatino}
\usepackage{latexsym}

\usepackage[mathscr]{eucal}
\usepackage{amsmath}
\usepackage{amsthm}
\usepackage{amsfonts}
\usepackage{amssymb}
\usepackage{amscd}
\usepackage{color}
\usepackage{graphicx}
\usepackage{graphics}
\usepackage{pifont}
\usepackage{subfigure}
\usepackage[makeroom]{cancel}
\usepackage[normalem]{ulem}
\usepackage[dvipsnames]{xcolor}
\usepackage{tikz}
\usepackage{bbm}

\theoremstyle{plain}
\newtheorem{theorem}{Theorem}

\newcommand{\ii}{\mathrm{i}}
\newcommand{\cH}{\mathcal{H}}

\newcommand{\ud}{\mathrm{d}}

\usepackage{newtxtext}       % 
\usepackage{newtxmath}  
\usepackage{tikz}

\title{Heat equation with inverse-square potential of bridging type across two half-lines}

\begin{document}
\author{Matteo Gallone\footnote{Department of Mathematics, Milan University, 
via Cesare Saldini 50, 20133 Milan (Italy) e-mail: matteo.gallone@unimi.it}  \and 
Alessandro Michelangeli\footnote{Institute for Applied Mathematics, and Hausdorff Center for Mathematics, University of Bonn  Endenicher Allee 60, D-53115 Bonn (Germany)
e-mail: michelangeli@iam.uni-bonn.de  }        \and
         Eugenio~Pozzoli\footnote{Institut de Math\'ematiques de Bourgogne, UMR 5584, CNRS, Universit\'e Bourgogne Franche-Comt\'e, F-21000 Dijon (France) e-mail: eugenio.pozzoli@u-bourgogne.fr
.}}
         
         \maketitle
\begin{abstract}
The heat equation with inverse square potential on both half-lines of $\mathbb{R}$ is discussed in the presence of \emph{bridging} boundary conditions at the origin. The problem is the lowest energy (zero-momentum) mode of the transmission of the heat flow across a Grushin-type cylinder, a generalisation of an almost Riemannian structure with compact singularity set. This and related models are reviewed, and the issue is posed of the analysis of the dispersive properties for the heat kernel generated by the underlying positive self-adjoint operator. Numerical integration is shown that provides a first insight and relevant qualitative features of the solution at later times.
\end{abstract}

\section{Introduction: the bridging heat equation in 1D.}\label{sec:intro}

% rigid \cyrzh\cyryo\cyrs\cyrt\cyrk\cyro\cyre  
% soft \cyrm\cyrya\cyrg\cyrk\cyro\cyre

For fixed $\alpha\in[0,1)$ we discuss in this note the following initial value problem in the unknown $u\equiv u(t,x)$, with $t\geqslant 0$ and $x\in\mathbb{R}\setminus\{0\}$:
\begin{equation}\label{eq:IVPgeneral}
\begin{cases}
 \;\displaystyle\frac{\partial u}{\partial t}-\frac{\partial^2 u}{\partial x^2}+\frac{\alpha(\alpha+2)}{4 x^2}\,u\;=\;0 \\
 \;u_0(t)^-\;=\;u_0(t)^+ & \!\!\!\!\!\!\!\!\!\!\!\!\!\!\!\!\!\!\!\!\!\!\!\!\!\!\!\!\textrm{where } u_0^\pm(t)\,:=\,\displaystyle\lim_{x\to 0^{\pm}}|x|^{\frac{\alpha}{2}}u(t,x)\,, \\
  \;u_1(t)^-\;=\;-u_1(t)^+ & \!\!\!\!\!\!\!\!\!\!\!\!\!\!\!\!\!\!\!\!\!\!\!\!\!\!\!\!\textrm{where } u_1^\pm(t)\,:=\,\displaystyle\lim_{x\to 0^{\pm}}|x|^{-(1+\frac{\alpha}{2})}(u(t,x)-|x|^{-\frac{\alpha}{2}}u_0^\pm(t))\,, \\
 \; u(0,x)\;=\;\varphi(x)& \!\!\!\!\!\!\!\!\!\!\!\!\!\!\!\!\!\!\!\!\!\!\!\!\!\!\!\!\textrm{where }\varphi\in L^2(\mathbb{R})\,,
\end{cases}
\end{equation}
seeking for solutions $u$ that, for (almost every) $t$, belong to $L^2(\mathbb{R})$. In the above formulation \eqref{eq:IVPgeneral} the existence of the limits $u_0^\pm$ and $u_1^\pm$ is part of the problem. We shall also consider the special case where the initial datum itself satisfies the very boundary conditions required at later times.

We are in particular concerned with the well-posedness of the problem and the dispersive properties of the solution(s).

At the same time, in this note we review the origin and meaning of the problem \eqref{eq:IVPgeneral} in the context of geometric quantum confinement or transmission across the metric's singularity for a particle constrained on a degenerate Riemannian manifold and only subject to the geometry of the constraining manifold, thus ``free to evolve'' over that manifold in analogy to a classical particle moving along geodesics.

The latter viewpoint is attracting an increasing amount of interest in recent years, making then natural to investigate the time-dependent equations arising in such contexts. Ours here is a `pilot' analysis of a more systematic study that unfolds ahead of us concerning dispersive and Strichartz estimates, and it has therefore the purpose of some propaganda and overview of the state of the art and on the future perspectives. Besides, here we only deal with the heat evolution, and not the Schr\"{o}dinger evolution, as we shall comment in due time.

 Prior to outlining the geometric background, let us comment on the structure of the problem \eqref{eq:IVPgeneral}. The considered PDE is a heat-type equation governed by the second order, elliptic (Schr\"{o}dinger) differential operator
 \begin{equation}\label{eq:diffop}
  -\frac{\ud^2}{\ud x^2}+\frac{C_\alpha}{\;x^2}\,,\qquad C_\alpha\,:=\,\frac{\alpha(\alpha+2)}{4 x^2}
 \end{equation}
 (the precise meaning of the parameter $\alpha$ and its presence through the coefficient $C_\alpha$ will be clear after discussing the parent geometric model). As such, the complete description of square-integrable solutions to the associated heat equation
 %, say, on the \emph{half-line} $(0,+\infty)$ 
 is achieved through a standard PDE analysis, once certain features of \eqref{eq:diffop} are known as a linear operator on $L^2(\mathbb{R})$.
 %$L^2(\mathbb{R}^+)$. 
 For concreteness, a limit-circle/limit-point argument \cite[Theorems X.11]{rs2} shows that when $C_\alpha
 \geqslant\frac{3}{4}$, i.e., $\alpha\in(-\infty,-3]\cup[1,+\infty)$, the linear operator \eqref{eq:diffop} minimally defined on smooth functions compactly supported away from $x=0$ is actually essentially self-adjoint on $L^2(\mathbb{R})$. Denoting its closure with $A$, one concludes that $A$ is a self-adjoint operator with strictly positive spectrum and domain $\mathcal{D}(A)$ that explicitly, when $C_\alpha>\frac{3}{4}$, is the Sobolev space $H^2_0(\mathbb{R})$. As a straightforward consequence of the abstract theory of differential equation on Hilbert space \cite[Proposition 6.6]{schmu_unbdd_sa}, one then concludes that the heat equation $\frac{\ud}{\ud t}u=-Au$ with initial datum $\varphi\in L^2(\mathbb{R})$ admits a unique solution in $C^1(\mathbb{R}_t^+,L^2(\mathbb{R}_x))$, with $u(t,\cdot)\in \mathcal{D}(A)$ at ever later $t>0$, explicitly given by $u(t,x)=(e^{-tA}\varphi)(x)$.

 In fact, it is worth recalling that the inverse square potential differential operator \eqref{eq:diffop} is greatly studied and deeply understood from many standpoints, in particular, both as far radial space-time (Strichartz) estimates are concerned both in the linear and non-linear Schr\"{o}dinger evolution (see, e.g., \cite{Burq-Planchon-Stalker-2003,M-2015-nonStrichartzHartree} and references therein), and as a Bessel operator on the $L^2$-space of the half-line (see, e.g., \cite{Derezinski-Georgescu-2021} and its precursors in that prolific research line).

 Yet, in addition to the differential side, the problem \eqref{eq:IVPgeneral} prescribes the solutions $u$ to satisfy certain boundary conditions at $x=0$. The first one, $g_0^+(t)=g_0^-(t)$ can be interpreted as the continuity of the function, up to the weight $|x|^{\frac{\alpha}{2}}$ that allows $u$ to have some degree of singularity at the origin; analogously, the condition $g_1^+(t)=-g_1^-(t)$ links the right and left derivative at zero, up to certain weights, and taken directionally from each side. In the regime $C_\alpha>\frac{3}{4}$ such conditions are obviously redundant, but when $C_\alpha\leqslant \frac{3}{4}$ an ad hoc analysis is needed to recognise that the prescribed behaviour at the origin expresses another condition of self-adjointness and positivity. Such an analysis has been carried out in several recent works \cite{Boscain-Prandi-JDE-2016,GMP-Grushin-2018,GMP-Grushin2-2020,PozzoliGru-2020volume,GM-Grushin3-2020} and is concisely reviewed in Section \ref{sec:Grucyl}. The net result, for the sake of the present discussion, is the following.
 
 \begin{theorem}[\cite{GMP-Grushin2-2020}]\label{th:bridging}
  Let $\alpha\in[0,1)$ and let $C_\alpha$ be given by \eqref{eq:diffop}.
  \begin{enumerate}
   \item[(i)] The space 
   \[
    \mathcal{D}\;:=\;\left\{ g\in L^2(\mathbb{R})\,\bigg|\,\left(-\frac{\ud^2}{\ud x^2}+\frac{C_\alpha}{\;x^2}\right)g\in L^2(\mathbb{R})\right\} 
   \]
  is a dense subspace of $L^2(\mathbb{R})$ and for every $g\in\mathcal{D}$ the following limits exist and are finite:
   \begin{equation}\label{eq:bilimitsg0g1}
  \begin{split}
   g_0^\pm\;&=\;\lim_{x\to 0^\pm} |x|^{\frac{\alpha}{2}}g^\pm(x)\,, \\
   g_1^\pm\;&=\;\lim_{x\to 0^\pm} |x|^{-(1+\frac{\alpha}{2})}\big(g^\pm(x)-g_0^\pm |x|^{-\frac{\alpha}{2}}\big)\,.
  \end{split}
 \end{equation}
 \item[(ii)] The operator
 \begin{equation}\label{eq:defineAbridge1d}
  \begin{split}
   \mathcal{D}\big(A_\alpha^{\mathrm{B}}\big)\;&=\;\{g\in\mathcal{D}\,|\, g_0^+=g_0^-\,,\; g_1^+=-g_1^-\}\,, \\
   A_\alpha^{\mathrm{B}} g\;&=\;- g''+ C_\alpha |x|^{-2} g
  \end{split}
 \end{equation}
  is self-adjoint on $L^2(\mathbb{R})$ and non-negative. Its spectrum is $[0,+\infty)$, and is all essential and absolutely continuous.
  \end{enumerate}
 \end{theorem}

 In Theorem \ref{th:bridging} we only consider the regime $\alpha\in[0,1)$. The remaining regime $\alpha\in[-3,0)$ is simply less relevant from the viewpoint of the underlying geometric model, as will be argued in Section \ref{sec:Grucyl}. And, as discussed above, when $\alpha\in(-\infty,-3)\cup(1,+\infty)$ one applies standard limit-point/limit-circle considerations.

  In view of Theorem \ref{th:bridging}, the initial value problem \eqref{eq:IVPgeneral} is immediately interpreted as the problem for the one-dimensional heat equation governed by the positive and self-adjoint operator $A_\alpha^{\mathrm{B}}$, and therefore it admits unique solution $u(t)=e^{-tA_\alpha^{\mathrm{B}}} \varphi$, again by abstract facts of differential equations on Hilbert space \cite[Proposition 6.6]{schmu_unbdd_sa}.

 The \emph{well-posedness} of \eqref{eq:IVPgeneral} is therefore fully controlled in $C^1(\mathbb{R}^+_t,L^2(\mathbb{R}_x))$ with $u(t,\cdot)\in\mathcal{D}(A_\alpha^{\mathrm{B}})$ at every $t>0$.

 The superscript `B' in $A_\alpha^{\mathrm{B}}$ is to refer to certain `\emph{bridging}' features of optimal transmission across the origin, allowing in a precise sense complete communication between the right and left half-line, induced by $A_\alpha^{\mathrm{B}}$, as compared to a whole family of similar transmission protocols.

 \begin{figure}[t!]
\begin{center}
\includegraphics[width=9cm]{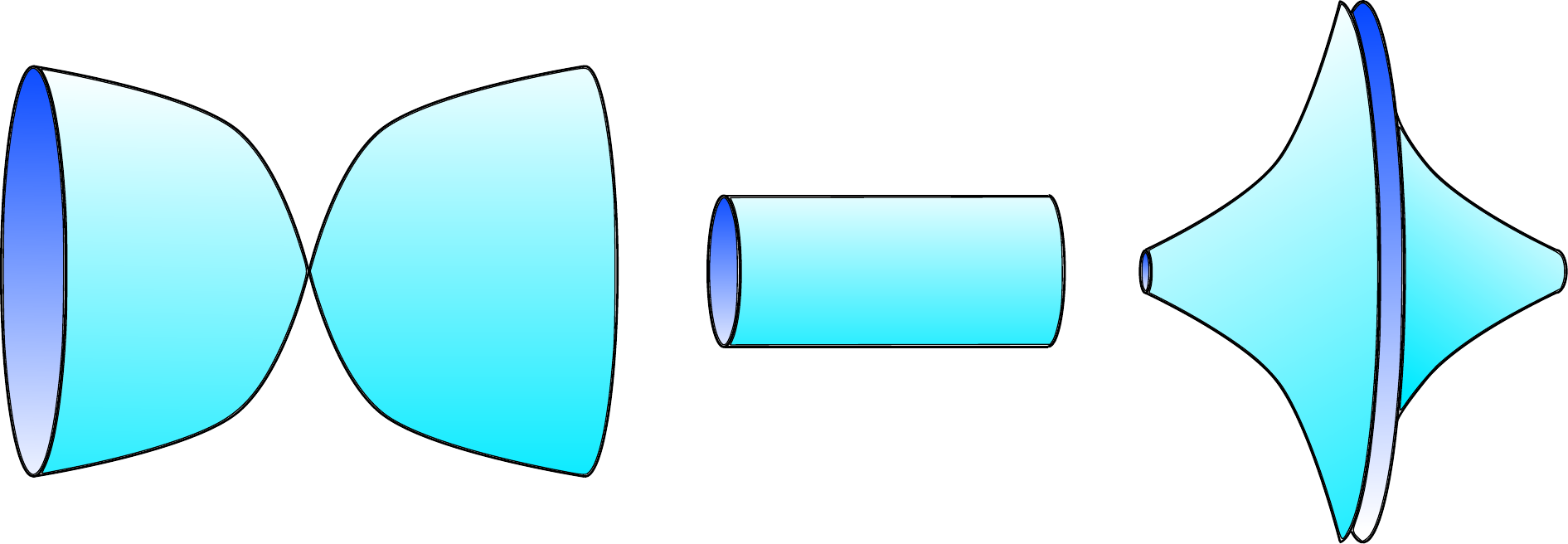}\caption{Manifold $M_\alpha$ for different $\alpha$: $\alpha<0$ (left), $\alpha=0$ (center) and $\alpha>0$ (right).}\label{fig:cylinders} 
\end{center}
\end{figure}

 Indeed, in Section \ref{sec:Grucyl} it will be recapped how initial value problems like \eqref{eq:IVPgeneral}, and their counterparts with the Schr\"{o}dinger equation, describe the heat-type or Schr\"{o}dinger-type propagation over a particular almost-Riemannian structure, customarily referred to as a `\emph{Grushin cylinder}', constituted by an infinite two-dimensional cylinder with a non-flat metric that becomes suitably singular on a given orthogonal section. Depending on the magnitude of the metric's singularity, which is quantified by the parameter $\alpha$, the transmission is either inhibited, so that a function initially supported on one half-cylinder remains confined in that half at later times, or on the opposite it is allowed, through a precise set of boundary conditions between the two halves. This can be qualitatively visualised as in Figure \ref{fig:cylinders} with cylinders that shrink to one point or get flattened in correspondence of the given singular section. \emph{One-dimensional} problems like \eqref{eq:IVPgeneral} emerge for the evolution on the lowest-energy mode, which corresponds to functions on the cylinder that are constant along the compact variable. It turns out that in certain regimes of metric's singularity (and $\alpha\in[0,1)$ is the physically most significant regime) an infinity of transmission protocols emerge, each characterised by suitable boundary conditions of self-adjointness, and each yielding a heat-type or Schr\"{o}dinger-type equation. Among them, the \emph{bridging} protocol described by the operator $A_\alpha^{\mathrm{B}}$ given by \eqref{eq:defineAbridge1d} displays distinguished features of optimal transmission. In practice -- see equations \eqref{eq:bridging_conditions}-\eqref{eq:bridging_transmission} below -- the heat (and Schr\"{o}dinger) equation of bridging type between the two half-lines is the one that describes a crossing at $x=0$ `without spatial filter' (continuity of the function) and `without energy filter' (the fraction of transmitted flux does not depend on the incident energy).

 The bridging protocol $A_\alpha^{\mathrm{B}}$ was first identified in the recent work \cite{Boscain-Prandi-JDE-2016}. The comparison analysis of the bridging protocol with respect to the whole family of the other physically meaningful ones was subsequently analysed in \cite{GMP-Grushin-2018,GM-Grushin3-2020}. Section \ref{sec:related-and-plane} reports on the recent literature of closely related models and results. Moreover, in \cite{Boscain-Prandi-JDE-2016,Boscain-Neel-2019} the bridging and some other protocols were analysed from the perspective of the conservation of the total heat, or equivalently, the perspective of infinite lifespan of the stochastic processes generated by such operators (stochastic completeness).

  It should be then sufficiently clear at this point that the initial value problem \eqref{eq:IVPgeneral} describes the low-energy transmission of bridging-type between the two halves of a Grushin cylinder with metric singularity at zero. In view of that, beside the already guaranteed well-posedness, \emph{dispersive properties} of the solution come to have great relevance in connection with the underlying physical transmission protocol, in particular $L^p$-$L^q$ estimates, smoothing estimates, and space-time (Strichartz) estimates for the heat semi-group associated with $A_\alpha^{\mathrm{B}}$.

  These are the analogue of the well-known estimates for the classical heat equation $(\frac{\partial}{\partial t}-\frac{\partial^2}{\partial x^2})u(t,x)=0$, $u(0,x)=\varphi(x)$, for which one has \cite[Section 2.2]{Wang-Zhaohui-Chengchun_HarmonicAnalysisI_2011}
 \begin{equation}\label{eq:heatdispest}
  \begin{split}
   \big\|\, u(t,\cdot)\,\big\|_{L^p}\;&\lesssim\;t^{-\frac{1}{2}(\frac{1}{r}-\frac{1}{p})}\|\varphi\|_{L^r}\qquad \quad 1\leqslant r\leqslant p\leqslant\infty\,, \\
   \Big\|\, \frac{\partial}{\partial x} u(t,\cdot)\,\Big\|_{L^p}\;&\lesssim\;t^{-\frac{1}{2}(\frac{1}{r}-\frac{1}{p}+1)}\|\varphi\|_{L^r}\qquad \, 1\leqslant r\leqslant p\leqslant\infty\,, \\
   \|\,u\,\|_{L^q(\mathbb{R}^+_t,L^p_x)}\;&\lesssim\;\|\varphi\|_{L^2}\qquad\qquad\qquad\;\;\, 2\leqslant p <\infty\,,\; q=\frac{4p}{p-2}\,, \\
   \Big\|\,  \frac{\partial}{\partial x}  u\,\Big\|_{L^2(\mathbb{R}^+_t,L^2_x)}\;&\lesssim\;\|\varphi\|_{L^2}\,.
  \end{split}
 \end{equation}

 Establishing the analogue of \eqref{eq:heatdispest} for the heat type semi-group $\exp(-t A_\alpha^{\mathrm{B}})$, and in fact eventually also for the Schr\"{o}dinger unitary group $\exp(-\ii t A_\alpha^{\mathrm{B}})$, as well as for the corresponding semi-groups and groups induced by other transmission protocols, and more generally on other geometries and classes of almost-Riemannian structures beyond the Grushin cylinder, appear to be one of most relevant challenges in this field, with abstract interest per se and impact on applications, of quantum control concern in the first place.

 In this respect we intend with this note to promote the above questions and advertise them for future investigations, in particular posing them in the rigorous context of geometric confinement and transmission protocols. Such an overview is given, as mentioned, in Sections \ref{sec:Grucyl} and \ref{sec:related-and-plane}.

 Last, in Section \ref{sec:num} we present a glance at numerical computations of the solution $u$ of \eqref{eq:IVPgeneral} when the initial datum is well localised on one half-line. The numerical evidence is strong on the dynamical formation of the bridging boundary conditions at $x=0$ at later times, and on a general behaviour that is qualitatively comparable to the classical heat propagation.

 In fact, at present no analytic computation is available of the heat (and, in the future, the Schr\"{o}dinger) propagator generated by $A_\alpha^{\mathrm{B}}$, and numerics is a first natural approach to infer meaningful properties to be rigorously proved in forthcoming investigations.
  
\section{A concise review of geometric confinement and transmission protocols in a Grushin cylinder}\label{sec:Grucyl}

 We have already anticipated that the problem \eqref{eq:IVPgeneral} provides the one-dimensional description for the heat flow across the singularity of a Grushin cylinder, and let us give in this Section a concise overview of the problem from that perspective.

 Grushin cylinders are Riemannian manifolds $M_\alpha\equiv(M,g_\alpha)$, with parameter $\alpha\in\mathbb{R}$, where
\begin{equation}
	M^\pm \; := \; \mathbb{R}^\pm_x \times \mathbb{S}^1_y \, \qquad \mathcal{Z} \;:=\; \{0\} \times \mathbb{S}^1_y \, , \qquad M\;:=\; M^+ \cup M^- \, 
\end{equation}
and with degenerate Riemannian metric 
\begin{equation}\label{eq:Grushin_Metric}
	g_\alpha \;:=\; \ud x \otimes \ud x +|x|^{-2\alpha} \ud y \otimes \ud y\,.
\end{equation}
Thus, $M_\alpha$ is a two-dimensional manifold built upon the cylinder $\mathbb{R} \times \mathbb{S}^1$, with singularity locus $\mathcal{Z}$ and incomplete Riemannian metric both on the right and the left half-cylinder $\mathbb{R}^\pm \times \mathbb{S}^1$ meaning that geodesics cross smoothly the singularity $\mathcal{Z}$ at finite times). The values $\alpha=-1$, $\alpha=0$, and $\alpha=1$ select, respectively, the flat cone, the Euclidean cylinder, and the standard `\emph{Grushin cylinder}' \cite[Chapter 11]{Calin-Chang-SubRiemannianGeometry}: in the latter case one has an `\emph{almost-Riemannian structure}'  on $\mathbb{R}\times\mathbb{S}^1=M^+\cup\mathcal{Z}\cup M^-$ in the rigorous sense of \cite[Sec.~1]{Agrachev-Boscain-Sigalotti-2008} or \cite[Sect.~7.1]{Prandi-Rizzi-Seri-2016}. Actually, $g_{\alpha}$ is defined as the unique metric for which the distribution of vector fields globally defined on $\mathbb{R}\times\mathbb{S}^1$ as
\begin{equation}\label{eq:frame}
X_1(x,y)\;:=\;\dfrac{\partial}{\partial x}\,, \qquad X^{(\alpha)}_2(x,y)\;:=\;|x|^\alpha \dfrac{\partial}{\partial y}
\end{equation}
is an orthonormal frame at every $(x,y)\in \mathbb{R}\times\mathbb{S}^1$: in this regard, the Grushin cylinder ($\alpha=1$) is a two-dimensional almost-Riemannian manifold of step two, meaning that 
$${\rm span}\Big\{X_1,X^{(1)}_2,\Big[X_1,X^{(1)}_2\Big]\Big\}\Big|_{(x,y)}\;=\;\mathbb{R}^2 \qquad \forall (x,y)\in  \mathbb{R}\times\mathbb{S}^1\,,$$
where $[X_1,X_2^\alpha]$ denotes the Lie brackets of vector fields. In fact, $M_\alpha$ is a hyperbolic manifold whenever $\alpha>0$, with Gaussian (sectional) curvature
\begin{equation}\label{eq:curvature}
 K_\alpha(x,y)\;=\;-\frac{\alpha(\alpha+1)}{x^2}\,.
\end{equation}

 To each $M_\alpha$ one naturally associates the Riemannian volume form
\begin{equation}\label{eq:volumeform}
 \mu_\alpha\;:=\;\mathrm{vol}_{g_\alpha}\;=\;\sqrt{\det g_\alpha}\,\ud x\wedge\ud y\;=\;|x|^{-\alpha}\,\ud x\wedge\ud y\,,
\end{equation}
 the Hilbert space 
 \begin{equation}\label{eq:Hspacealpha}
 \cH_\alpha\;:=\;L^2(M,\ud\mu_\alpha)\,,
\end{equation}
understood as the completion of $C^\infty_c(M)$ with respect to the scalar product
\begin{equation}
 \langle \psi,\varphi\rangle_{\alpha}\;:=\;\iint_{\mathbb{R}\times\mathbb{S}^1}\overline{\psi(x,y)}\,\varphi(x,y)\,\frac{1}{|x|^{\alpha}}\,\ud x\,\ud y\,,
\end{equation}
 and the (Riemannian) Laplace-Beltrami operator $\Delta_{\mu_\alpha}:=\mathrm{div}_{\mu_\alpha}\circ\nabla$ acing on functions over $M_\alpha$. A standard computation (see, e.g. \cite[Sect.~2]{GMP-Grushin-2018}) yields explicitly
\begin{equation}\label{eq:IntroLaplaceBeltrami}
	\Delta_{\mu_\alpha} \; = \; \frac{\partial^2}{\partial x^2} +|x|^{2\alpha}\frac{\partial^2}{\partial y^2} -\frac{\alpha}{|x|} \frac{\partial}{\partial x}\,.
\end{equation}

 As a linear operator on $\cH_\alpha$ one minimally defines $\Delta_{\mu_\alpha}$ on the dense subspace of $L^2(M,\ud\mu_\alpha)$-functions that are smooth and compactly supported away from the metric's singularity locus $\mathcal{Z}$, thus introducing
\begin{equation}\label{Halpha}
	H_\alpha\;:=\; - \Delta_{\mu_\alpha} \, , \qquad \mathcal{D}(H_\alpha) \;:= \; C^\infty_c(M)\,.
\end{equation}
 The Green identity implies that $H_\alpha$ is symmetric and non-negative. One has the following classification.

 \begin{theorem}[\cite{Boscain-Laurent-2013,Boscain-Prandi-JDE-2016,GMP-Grushin-2018}]\label{thm:cylinder}
  \begin{enumerate}
   \item[(i)] If $\alpha\in(-\infty,-3]\cup[1,+\infty)$, then the operator $H_\alpha$ is essentially self-adjoint.
   \item[(ii)] If $\alpha\in(-3,-1]$, then $H_\alpha$ is not essentially self-adjoint with deficiency index $2$.
   \item[(iii)] If $\alpha\in(-1,+1)$, then $H_\alpha$ is not essentially self-adjoint and it has infinite deficiency index.
  \end{enumerate}
  \end{theorem}
  The same holds, separately, for the symmetric operators $H^\pm_\alpha$ minimally defined on the $L^2$-space of each half-cylinder.

  With respect to the Hilbert space orthogonal decomposition
  \begin{equation}
   \cH_\alpha\;=\;L^2(M,\ud\mu_\alpha)\;\cong\;L^2(M^-,\ud\mu_\alpha\oplus L^2(M^+,\ud\mu_\alpha)
  \end{equation}
  the operator $H_\alpha$ is reduced as $H_\alpha=H_\alpha^-\oplus H_\alpha^+$, and therefore in the regime of essential self-adjointness its closure is the reduced, non-negative, self-adjoint operator $\overline{H_\alpha}=\overline{H_\alpha^-}\oplus \overline{H_\alpha^+}$. This implies that both the Schr\"{o}dinger equation $\partial_t u=-\ii\overline{H_\alpha} u$ and the heat equation $\partial_tu=-\overline{H_\alpha} u$ decompose to uncoupled equations on each half-cylinder, or, better to say, group and semi-group decompose, respectively, as $e^{-\ii t \overline{H_\alpha}}=e^{-\ii t \overline{H_\alpha^-}}\oplus e^{-\ii t \overline{H_\alpha^+}}$ and $e^{-t \overline{H_\alpha}}=e^{- t \overline{H_\alpha^-}}\oplus e^{- t \overline{H_\alpha^+}}$, with the consequence that an initial datum supported, say, only on one half, keeps evolving in that half at each later time. This phenomenon is customarily referred to as `\emph{heat-geometric confinement}' and `\emph{quantum-geometric confinement}', respectively, to emphasise the sole effect of the geometry (meaning that we are not considering any potential energy on the manifold, but only the kinetic one), with no coupling boundary conditions -- hence no interaction -- declared at $\mathcal{Z}$. Quantum-mechanically, in this regime of the Grushin metric, a quantum particle constrained on $M_\alpha$ and left `free' to evolve only under the effect of the underlying geometry, never happen to cross the singularity locus $\mathcal{Z}$.

  The scenario becomes much more diversified when $H_\alpha$ is \emph{not} essentially self-adjoint and therefore admits non-trivial self-adjoint extensions. Our regime of interest includes $\alpha\in(0,1)$, the sub-case of greatest relevance because it corresponds to an actual local singularity (and not vanishing) of the metric $g_\alpha$, and for the purposes of this note we shall only consider such $\alpha$'s. Qualitatively analogous results can be established in the remaining non-self-adjoint regime $\alpha\in(-3,0)$.

  There is in fact a giant family of inequivalent self-adjoint realisations of $H_\alpha$ when $\alpha\in[0,1)$, as the deficiency index is infinite. Each one is characterised by boundary conditions at $\mathcal{Z}$ that prescribe a one-sided or two-sided interaction with the boundary, or more generally a protocol of left$\leftrightarrow$right transmission. Such a family include physically unstable realisations (those that are not lower semi-bounded), as well as a huge amount of unphysical realisations, such as those with non-local boundary conditions at $\mathcal{Z}$.

  An extensive and fairly explicit classification of physical extensions of $H_\alpha$ was recently completed in \cite{GMP-Grushin2-2020}. 
  %Here and in the following we 

  \begin{theorem}[\cite{GMP-Grushin2-2020}]\label{thm:H_alpha_fibred_extensions}
 Let $\alpha\in[0,1)$. $H_\alpha$ defined in \eqref{Halpha} admits, among others, the following families of self-adjoint extensions with respect to $L^2(M,\ud\mu_\alpha)$:
 \begin{itemize}
  \item \underline{Friedrichs extension}: $H_{\alpha,\mathrm{F}}$;
  \item \underline{Family $\mathrm{I_R}$}: $\{H_{\alpha,\mathrm{R}}^{[\gamma]}\,|\,\gamma\in\mathbb{R}\}$;
  \item \underline{Family $\mathrm{I_L}$}: $\{H_{\alpha,\mathrm{L}}^{[\gamma]}\,|\,\gamma\in\mathbb{R}\}$;
  \item \underline{Family $\mathrm{II}_a$} with $a\in\mathbb{C}$: $\{H_{\alpha,a}^{[\gamma]}\,|\,\gamma\in\mathbb{R}\}$;
  \item \underline{Family $\mathrm{III}$}: $\{H_{\alpha}^{[\Gamma]}\,|\,\Gamma\equiv(\gamma_1,\gamma_2,\gamma_3,\gamma_4)\in\mathbb{R}^4\}$.
 \end{itemize}
 Each member of any such family acts precisely as the differential operator $-\Delta_{\mu_\alpha}$ on a domain of functions $f\in L^2(M,\ud\mu_\alpha)$ satisfying the following properties.
  \begin{enumerate}
  \item[(i)] \underline{Integrability and regularity}:
  \begin{equation}\label{eq:DHalpha_cond1}
  \sum_{\pm}\;\iint_{\mathbb{R}_x^\pm\times\mathbb{S}^1_y}\big|(\Delta_{\mu_\alpha}f)(x,y)\big|^2\,\ud\mu_\alpha(x,y)\;<\;+\infty\,.
 \end{equation}
  \item[(ii)] \underline{Boundary condition}: The limits
 \begin{eqnarray}
  f_0^\pm(y)&=&\lim_{x\to 0^\pm}f(x,y) \label{eq:DHalpha_cond2_limits-1}\\
  f_1^\pm(y)&=&\pm(1+\alpha)^{-1}\lim_{x\to 0^\pm}\Big(\frac{1}{\:|x|^\alpha}\,\frac{\partial f(x,y)}{\partial x}\Big) \label{eq:DHalpha_cond2_limits-2}
  \end{eqnarray}
 exist and are finite for almost every $y\in\mathbb{S}^1$, and depending on the considered type of extension, and for almost  every $y\in\mathbb{R}$,
 \begin{eqnarray}
  f_0^\pm(y)\,=\,0 \qquad \quad\;\;& & \textrm{if }\;  f\in\mathcal{D}(H_{\alpha,\mathrm{F}})\,, \label{eq:DHalpha_cond3_Friedrichs}\\
  \begin{cases}
   \;f_0^-(y)= 0  \\
   \;f_1^+(y)=\gamma f_0^+(y)
  \end{cases} & & \textrm{if }\;  f\in\mathcal{D}(H_{\alpha,\mathrm{R}}^{[\gamma]})\,, \\
   \begin{cases}
   \;f_1^-(y)=\gamma f_0^-(y) \\
   \;f_0^+(y)= 0 
  \end{cases} & & \textrm{if }\;  f\in\mathcal{D}(H_{\alpha,\mathrm{L}}^{[\gamma]}) \,, \label{eq:DHalpha_cond3_L}\\
     \begin{cases}
   \;f_0^+(y)=a\,f_0^-(y) \\
   \;f_1^-(y)+\overline{a}\,f_1^+(y)=\gamma f_0^-(y)
  \end{cases} & & \textrm{if }\;  f\in\mathcal{D}(H_{\alpha,a}^{[\gamma]})\,, \label{eq:DHalpha_cond3_IIa} \\
   \begin{cases}
   \;f_1^-(y)=\gamma_1 f_0^-(y)+(\gamma_2+\ii\gamma_3) f_0^+(y) \\
   \;f_1^+(y)=(\gamma_2-\ii\gamma_3) f_0^-(y)+\gamma_4 f_0^+(y)
  \end{cases} & & \textrm{if }\;  f\in\mathcal{D}(H_{\alpha}^{[\Gamma]})\,. \label{eq:DHalpha_cond3_III}
 \end{eqnarray} 
 \end{enumerate} 
%      Moreover,
%  \begin{equation}\label{eq:traceregularity}
%   f_0^\pm \in H^{s_{0,\pm}}(\mathbb{S}^1, \ud y)\qquad\textrm{ and }\qquad f_1^\pm\in H^{s_{1,\pm}}(\mathbb{S}^1,\ud y)
%  \end{equation}
%  with
%  \begin{itemize}
%  	\item $s_{1,\pm}=\frac{1}{2}\frac{1-\alpha}{1+\alpha}$\qquad\qquad\qquad\qquad\qquad\; for the Friedrichs extension,
%  	\item $s_{1,-}=\frac{1}{2}\frac{1-\alpha}{1+\alpha}$, $s_{0,+}=s_{1,+}=\frac{1}{2}\frac{3+\alpha}{1+\alpha}$ \quad\quad for extensions of type $\mathrm{I_R}$\,,
%  	\item  $s_{1,+}=\frac{1}{2}\frac{1-\alpha}{1+\alpha}$, $s_{0,-}=s_{1,-}=\frac{1}{2}\frac{3+\alpha}{1+\alpha}$ \quad\quad for extensions of type $\mathrm{I_L}$\,,
%  	\item $s_{1,\pm}=s_{0,\pm}=\frac{1}{2}\frac{1-\alpha}{1+\alpha}$ \qquad\qquad\qquad \;\;\;\; for extensions of type $\mathrm{II}_a$\,,
%  	\item $s_{1,\pm}=s_{0,\pm}=\frac{1}{2}\frac{3+\alpha}{1+\alpha}$ \qquad\qquad\qquad \;\;\;\, for extensions of type $\mathrm{III}$\,.
%  \end{itemize}
\end{theorem}

  One can further select those extensions that are non-negative and then induce a heat type flow.
  
 \begin{theorem}[\cite{GM-Grushin3-2020}]\label{thm:positivity}
\begin{itemize}
		\item The Friedrichs extension $H_{\alpha,\mathrm{F}}$ is non-negative.
		\item Extensions in the family $\mathrm{I_R}$, $\mathrm{I_L}$, and $\mathrm{II}_a$, $a\in\mathbb{C}$, are non-negative if and only if $\gamma\geqslant 0$.
		\item Extensions in the family $\mathrm{III}$ are non-negative if and only if so is the matrix
		\[
		 \widetilde{\Gamma}\;:=\;\begin{pmatrix}
		  \gamma_1 & \gamma_2+\ii\gamma_3 \\
		  \gamma_2-\ii\gamma_3 & \gamma_4
		 \end{pmatrix},
		\]
                i.e., if and only if $\gamma_1 + \gamma_4 > 0$ and $\gamma_1 \gamma_4\geqslant \gamma_2^2 + \gamma_3^2$.
\end{itemize}
\end{theorem}

 A customary \emph{quantum-mechanical} quantification of the transmission modelled by each extension is the fraction of Schr\"{o}dinger flux that gets transmitted vs reflected when a beam of particles are shot free from infinity towards $\mathcal{Z}$. This analysis, albeit in a \emph{Schr\"{o}dinger equation} framework, elucidates the qualitative properties of the crossing at $x=0$ and was recently done in \cite{GM-Grushin3-2020}. Intuitively speaking, far away from $\mathcal{Z}$ the metric tends to become flat and the action $-\Delta_{\mu_\alpha}$ of each self-adjoint `free Hamiltonian' tends to resemble that of the free Laplacian $-\Delta$, plus the correction due to the $(|x|^{-1}\partial_x)$-term, on wave functions $f(x,y)$ that are constant in $y$. This suggests that at very large distances a quantum particle evolves free from the effects of the underlying geometry, and one can speak of scattering states of energy $E>0$. The precise shape of the wave function $f_{\mathrm{scatt}}$ of such a scattering state can be easily guessed to be of the form
\begin{equation}\label{eq:fscatt}
 f_{\mathrm{scatt}}(x,y)\;\sim\;|x|^{\frac{\alpha}{2}}e^{\pm\ii x\sqrt{E}}\qquad \textrm{as }|x|\to +\infty\,.
\end{equation}
Indeed, $-\Delta_{\mu_\alpha} f_{\mathrm{scatt}}\sim Ef_{\mathrm{scatt}}+\frac{\alpha(2+\alpha)}{4|x|^2}f_{\mathrm{scatt}}$, that is, up to a very small $O(|x|^{-2})$-correction, $f_{\mathrm{scatt}}$ is a generalised eigenfunction of $-\Delta_{\mu_\alpha}$ with eigenvalue $E$. All this can be fully justified on rigorous grounds \cite{GM-Grushin3-2020} and leads naturally to the definition of the `\emph{transmission coefficient}' and `\emph{reflection coefficient}' for the scattering, namely the spatial density of the transmitted flux and the reflected flux, normalised with respect to the density of the incident flux. Obviously, no scattering across the singularity occurs for Friedrichs, or type-$\mathrm{I}_\mathrm{R}$, or type-$\mathrm{I}_\mathrm{L}$ quantum protocols, whereas in type-$\mathrm{II}_a$ scattering one obtains the following (analogous conclusions can be made for type-$\mathrm{III}$ scattering).

 \begin{theorem}[\cite{GM-Grushin3-2020}]\label{thm:scattering} Let $\alpha\in[0,1)$, $a\in\mathbb{C}$, $\gamma\in\mathbb{R}$.
 The transmission coefficient $T_{\alpha,a,\gamma}(E)$ and the reflection coefficient $R_{\alpha,a,\gamma}(E)$ at given energy $E>0$ for the Schr\"{o}dinger transmission protocol governed by $H_{\alpha,a}^{[\gamma]}$ are given by
  \begin{equation}\label{eq:TR-IIa-thm}
  \begin{split}
   T_{\alpha,a,\gamma}(E)\;&=\;\left|\frac{\,E^{\frac{1+\alpha}{2}}(1+e^{\ii\pi\alpha})\,\Gamma(\frac{1-\alpha}{2})\,\overline{a}}{\,E^{\frac{1+\alpha}{2}}\Gamma(\frac{1-\alpha}{2})(1+|a|^2)+\ii\,\gamma\,2^{1+\alpha}e^{\ii\frac{\pi}{2}\alpha}\Gamma(\frac{3+\alpha}{2})} \right|^2 \,, \\
   R_{\alpha,a,\gamma}(E)\;&=\;\left|\frac{\,E^{\frac{1+\alpha}{2}}\Gamma(\frac{1-\alpha}{2})\,(1-|a|^2\,e^{\ii\pi\alpha})+\ii\,\gamma\,2^{1+\alpha}e^{\ii\frac{\pi}{2}\alpha}\Gamma(\frac{3+\alpha}{2})}{\,E^{\frac{1+\alpha}{2}}\Gamma(\frac{1-\alpha}{2})(1+|a|^2)+\ii\,\gamma\,2^{1+\alpha}e^{\ii\frac{\pi}{2}\alpha}\Gamma(\frac{3+\alpha}{2})}\right|^2\,.
  \end{split}
 \end{equation}
  They satisfy 
  \begin{equation}\label{eq:TpEi1-thm}
  T_{\alpha,a,\gamma}(E)+R_{\alpha,a,\gamma}(E)\;=\;1\,,
 \end{equation}
 and when $\gamma=0$ they are independent of $E$ .
 The scattering is reflection-less ($R_{\alpha,a,\gamma}(E)=0$) when
 \begin{equation}\label{eq:Etransm-thm}
   E\;=\;\bigg(\frac{\,2^{1+\alpha}\,\gamma\,\Gamma(\frac{3+\alpha}{2})\sin\frac{\pi}{2}\alpha\,}{\,\Gamma(\frac{1-\alpha}{2})(1-\cos\pi\alpha)\,}\bigg)^{\!\frac{2}{1+\alpha}},	
 \end{equation}
 provided that $\alpha\in(0,1)$, $|a|=1$, and $\gamma>0$. In the high energy limit the scattering is independent of the extension parameter $\gamma$ and one has
 \begin{equation}\label{eq:highenergyscatt-thm}
  \begin{split}
   \lim_{E\to+\infty}T_{\alpha,a,\gamma}(E)\;&=\;\frac{\,2\,|a|^2(1+\cos\pi\alpha)}{(1+|a|^2)^2}\,, \\
   \lim_{E\to+\infty}R_{\alpha,a,\gamma}(E)\;&=\;\frac{\,1+|a|^4-2|a|^2\cos\pi\alpha}{(1+|a|^2)^2}\,,
  \end{split}
 \end{equation}
 whereas in the low energy limit, for $\gamma\neq 0$,
  \begin{equation}\label{eq:lowenergyscatt-thm}
  \begin{split}
   \lim_{E\downarrow 0}T_{\alpha,a,\gamma}(E)\;&=\;0\,, \\
   \lim_{E\downarrow 0}R_{\alpha,a,\gamma}(E)\;&=\;1\,.
  \end{split}
 \end{equation}
\end{theorem}

 Upon inspection of the boundary conditions \eqref{eq:DHalpha_cond3_Friedrichs}-\eqref{eq:DHalpha_cond3_III} one sees that the type-$\mathrm{II}_a$ extension $H_{\alpha,a}^{[\gamma]}$ with $a=1$ and $\gamma=0$ imposes the local behaviour
\begin{equation}\label{eq:bridging_conditions}
 \begin{split}
  \lim_{x\to 0^-}f(x,y)\;&=\;\lim_{x\to 0^+}f(x,y) \\
  \lim_{x\to 0^-}\Big(\frac{1}{\:|x|^\alpha}\,\frac{\partial f(x,y)}{\partial x}\Big)\;&=\;\lim_{x\to 0^+}\Big(\frac{1}{\:|x|^\alpha}\,\frac{\partial f(x,y)}{\partial x}\Big)\,,
 \end{split}
\end{equation}
 namely the distinguished feature of having a domain of functions that are continuous across the Grushin singularity, together with their weighted derivative. $H_{\alpha,1}^{[0]}$ is called the `\emph{bridging extension}' of $H_\alpha$, and for it we shall simply write $H_\alpha^{\mathrm{B}}$. In view of the results reviewed so far, the transmission modelled by the bridging extension
 \begin{itemize}
 \item has \emph{no spatial filter} in the sense of \eqref{eq:bridging_conditions} (in fact, all type-$\mathrm{II}_a$ protocols with $a=1$ impose local continuity at $\mathcal{Z}$; quantum-mechanically this is interpreted as a lack of jump in the particle's probability density from one side to the other of the singularity),
 \item and has \emph{no energy filter} in the Schr\"{o}dinger scattering, indeed, $H_\alpha^{\mathrm{B}}$ and all type-$\mathrm{II}_a$ protocols with $\gamma=0$ induce a scattering where the fraction of transmitted and reflected flux does not depend on the incident energy (see \eqref{eq:TR-IIa-thm} above),  
 \begin{equation}\label{eq:bridging_transmission}
 \begin{split}
   T_\alpha^{\mathrm{B}}\;:=\;T_{\alpha,1,0}(E)\;&=\;\frac{1}{2}\,(1+\cos\pi\alpha)\,, \\
   R_\alpha^{\mathrm{B}}\;:=\; R_{\alpha,1,0}(E)\;&=\;\frac{1}{2}\,(1-\cos\pi\alpha)\,,
  \end{split}
\end{equation}
 meaning that the singularity does not act as a filter in the energy.  
\end{itemize}

 The overall picture surveyed so far poses naturally the question of the analysis of the heat type flow generated by the \emph{positive} and \emph{self-adjoint} realisations of the Laplace-Beltrami operator on $\cH_\alpha$ (Theorem \ref{thm:positivity}), as well as the Schr\"{o}dinger type flow generated by \emph{self-adjoint realisations} (Theorem \ref{thm:H_alpha_fibred_extensions}), let alone the study of \emph{non-linear} heat and Schr\"{o}dinger equations on $M_\alpha$ with linear part given by a self-adjoint Laplace-Beltrami operator. This appears to be a completely uncharted territory of notable relevance in abstract terms and for applications. The gap between such future goals and the current knowledge is a lack of informative characterisation of the heat and Schr\"{o}dinger propagator's kernel.

 To complete the present review, let us make the connection explicit between the two-dimensional heat type equation induced by $H_\alpha^{\mathrm{B}}$ and the one-dimensional problem \eqref{eq:IVPgeneral}.

 This is done \cite{Boscain-Prandi-JDE-2016,GMP-Grushin2-2020} by means of the canonical Hilbert space unitary isomorphism $\cH_\alpha\xrightarrow{\cong}\cH$, where
\begin{equation}\label{eq:Hxispace}
 \begin{split}
   \cH\;&:=\;\mathcal{F}_2U_\alpha L^2(M,\ud\mu_\alpha)\;\cong\;\ell^2(\mathbb{Z},L^2(\mathbb{R},\ud x))\;\cong\;\cH^-\oplus\cH^+\;\cong\;\bigoplus_{k\in\mathbb{Z}}\;\mathfrak{h}\,, \\
    \mathfrak{h}\;&:=\;L^2(\mathbb{R}^-,\ud x)\oplus L^2(\mathbb{R}^+,\ud x)\;\cong\;L^2(\mathbb{R},\ud x)\,,
 \end{split}
\end{equation}
 recalling that 
 \begin{equation}
  \cH_\alpha\;\cong\;L^2(M^-,\ud\mu_\alpha)\oplus L^2(M^+,\ud\mu_\alpha)\,,
 \end{equation}
 and where the unitary transformations $U_\alpha:=U_\alpha^-\oplus U_\alpha^+$ and $\mathcal{F}_2:=\mathcal{F}_2^-\oplus\mathcal{F}_2^+$ are defined, respectively, as
\begin{equation}\label{eq:unit1}
\begin{split}
 U_\alpha^\pm:L^2(\mathbb{R}^\pm\times\mathbb{S}^1,|x|^{-\alpha}\ud x\ud y)&\stackrel{\cong}{\longrightarrow}L^2(\mathbb{R}^\pm\times\mathbb{S}^1,\ud x\ud y) \\
 f &\; \mapsto\;\phi\;:=\;  |x|^{-\frac{\alpha}{2}}f\,,
\end{split}
\end{equation}
and
\begin{equation}\label{eq:defF2}
 \begin{split}
  \mathcal{F}_2^{\pm}:L^2(\mathbb{R}^\pm\times\mathbb{S}^1,\ud x\ud y)&\stackrel{\cong}{\longrightarrow}
 L^2(\mathbb{R}^\pm,\ud x)\otimes\ell^2(\mathbb{Z})\,, \\
  \phi &\;\mapsto\;\psi\;\equiv\;(\psi_k)_{k\in\mathbb{Z}}\,, \\
  e_k(y)\;:=\;\frac{e^{\ii k y}}{\sqrt{2\pi}}\,,&\qquad \psi_k(x)\,:=\int_0^{2\pi}\overline{e_k(y)}\,\phi(x,y)\,\ud y\,,\qquad x\in\mathbb{R}^{\pm}
 \end{split}
\end{equation}
(thus, $\phi(x,y)=\sum_{k\in\mathbb{Z}}\psi_k(x)e_k(y)$ in the $L^2$-convergent sense). This provides, up to isomorphism, the orthogonal sum decomposition of the Hilbert space of interest into identical `bilateral' \emph{fibres} $\mathfrak{h}=L^2(\mathbb{R}^-,\ud x)\oplus L^2(\mathbb{R}^+,\ud x)\cong L^2(\mathbb{R},\ud x)$. The decomposition is discrete, as a consequence of having taken the Fourier transform $\mathcal{F}_2$ only in the compact variable $y$.

 \begin{theorem}[\cite{GMP-Grushin2-2020}]\label{thm:connectionB}
  Let $\alpha\in[0,1)$. Through the isomorphism \eqref{eq:Hxispace} the self-adjoint bridging operator $H_\alpha^{\mathrm{B}}$ on $\cH_\alpha=L^2(M,\ud\mu_\alpha)$ is unitarily equivalent to the self-adjoint operator $\mathscr{H}_{\alpha}^{\mathrm{B}}$ on $\cH\cong\ell^2(\mathbb{Z},L^2(\mathbb{R}))$, namely
    \begin{equation}\label{eq:unitequiv1}
   H_\alpha^{\mathrm{B}}\;=\;(U_\alpha)^{-1}(\mathcal{F}_2)^{-1}\mathscr{H}_{\alpha}^{\mathrm{B}}\,\mathcal{F}_2\,U_\alpha\,,
  \end{equation}
  where
  \begin{equation}\label{eq:unitequiv2}
   \mathscr{H}_{\alpha}^{\mathrm{B}}\;=\;\bigoplus_{k\in\mathbb{Z}} A_{\alpha}(k)
  \end{equation}
  and each $A_\alpha(k)$ is the self-adjoint operator on $L^2(\mathbb{R})$ given by
  \begin{equation}\label{eq:Afibre}
 \begin{split}
  \mathcal{D}(A_{\alpha}(k))\;&=\;
  \left\{\!\!
  \begin{array}{c}
   g=g^-\oplus g^+\,,\; g^\pm\in L^2(\mathbb{R}^\pm,\ud x)\;\;\textrm{such that} \\
   \big(-\frac{\ud^2}{\ud x^2}+k^2 |x|^{2\alpha}+\frac{\,\alpha(2+\alpha)\,}{4x^2}\big)g^\pm\in L^2(\mathbb{R}^\pm,\ud x) \\
   g_0^-=g_0^+\,,\quad g_1^-=-g_1^+
  \end{array}
  \!\!\right\}, \\
   A_{\alpha}(k)g\;&=\;\bigoplus_{\pm}\Big(-\frac{\ud^2}{\ud x^2}+k^2 |x|^{2\alpha}+\frac{\,\alpha(2+\alpha)\,}{4x^2}\Big) g^\pm\,,
 \end{split}
\end{equation}
 where $g_0^\pm,g_1^\pm\in\mathbb{C}$ are the existing and finite limits
 \begin{equation}\label{eq:bilimitsg0g1-gen}
  \begin{split}
   g_0^\pm\;&=\;\lim_{x\to 0^\pm} |x|^{\frac{\alpha}{2}}g(x) \\
   g_1^\pm\;&=\;\lim_{x\to 0^\pm} |x|^{-(1+\frac{\alpha}{2})}\big(g(x)-g_0^\pm |x|^{-\frac{\alpha}{2}}\big)\,.
  \end{split}
 \end{equation}
 \end{theorem}

 In Theorem \ref{thm:connectionB} the existence and finiteness of the limits \eqref{eq:bilimitsg0g1-gen} is guaranteed by the distributional constraint $\big(-\frac{\ud^2}{\ud x^2}+k^2 |x|^{2\alpha}+\frac{\,\alpha(2+\alpha)\,}{4x^2}\big)g^\pm\in L^2(\mathbb{R}^\pm,\ud x)$. A completely analogous unitary equivalence and fibred decomposition like \eqref{eq:unitequiv1}-\eqref{eq:unitequiv2} holds for all other self-adjoint realisations of the Laplace-Beltrami operator on Grushin cylinder, as classified in Theorem \ref{thm:H_alpha_fibred_extensions} \cite{GMP-Grushin2-2020}.
 
 Each $A_\alpha(k)$ is the $k$-th transversal momentum mode of the operator $H_\alpha^{\mathrm{B}}$ on cylinder, in the sense of the isomorphism \eqref{eq:Hxispace}, namely with respect to the momentum conjugate to the $y$-variable. By compactness, these are discrete modes and, as seen from \eqref{eq:Afibre}, the boundary condition at $x=0$ has the \emph{same} form ($g_0^-=g_0^+$, $g_1^-=-g_1^+$) in \emph{each} mode, and moreover it \emph{does not couple distinct modes}. Because of this structure, the bridging operator $H_\alpha^{\mathrm{B}}$ is said to be `\emph{uniformly fibred}', and in fact all other extensions classified in Theorem \ref{thm:H_alpha_fibred_extensions} are uniformly fibred too \cite{GMP-Grushin2-2020}. Uniformly fibred extensions generate a heat or Schr\"{o}dinger flow that is reduced into the discrete modes $k$.

 A careful spectral analysis \cite{GM-Grushin3-2020} shows that for each (uniformly fibred) extension from Theorem \ref{thm:H_alpha_fibred_extensions}, the transversal momentum modes are energetically increasingly ordered in the sense of increasing $|k|$, meaning in particular that the zero-th mode is the lowest energy one, and that for the bridging operator all modes have only non-negative, essential, absolutely continuous spectrum.

 Comparing \eqref{eq:Afibre} with \eqref{eq:defineAbridge1d} one recognises that $A_\alpha(0)=A_\alpha^{\mathrm{B}}$. This and the considerations made in Section \ref{sec:intro} finally show that the heat flow generated by the bridging operator $H_\alpha^{\mathrm{B}}$ starting with a function $f_{\mathrm{iniz}}$ on the cylinder which belongs to the zero-th transversal momentum mode and therefore is constant in $y$, say, $f_{\mathrm{iniz}}(x,y)=\varphi(x)$, produces at times $t>0$ and evoluted function
 \begin{equation}
  f(t;x,y)\;=\;u(t,x)
 \end{equation}
 (still belonging to the zero mode) where $u$ solves the one-dimensional initial value problem \eqref{eq:IVPgeneral} with initial datum $\varphi$.

\section{Related settings: Grushin planes and almost Riemannian manifolds}\label{sec:related-and-plane}

  The subject of geometric quantum confinement away from the metric's singularity, and transmission across it, for quantum particles or for the heat flow on degenerate Riemannian manifolds is experiencing a fast growth in the recent years. Such themes are particularly active with reference to Grushin structures on cylinder, cone, and plane \cite{Boscain-Laurent-2013,Boscain-Prandi-Seri-2014-CPDE2016,Boscain-Prandi-JDE-2016,GMP-Grushin-2018,PozzoliGru-2020volume,Boscain-Neel-2019,Boscain-Beschastnnyi-Pozzoli-2020,IB-2021}, as well as, more generally, on two-dimensional orientable compact almost-Riemannian manifolds of step two \cite{Boscain-Laurent-2013}, $d$-dimensional regular almost-Riemannian and sub-Riemannian manifolds \cite{Prandi-Rizzi-Seri-2016,Franceschi-Prandi-Rizzi-2017}.

  Of significant relevance is the counterpart model to the Grushin-type cylinder, but in the lack of compact variable. This leads to related almost Riemannian structures called `\emph{Grushin-type planes}'. In complete analogy to Section \ref{sec:Grucyl}, these are Riemannian manifolds $M_\alpha\equiv(M,g_\alpha)$, for some $\alpha\in\mathbb{R}$, where now
  \begin{equation}
	M^\pm \; := \; \mathbb{R}^\pm_x \times \mathbb{R}_y \, \qquad \mathcal{Z} \;:=\; \{0\} \times \mathbb{R}_y \, , \qquad M\;:=\; M^+ \cup M^- \, 
\end{equation}
and again with degenerate Riemannian metric 
\begin{equation}
	g_\alpha \;:=\; \ud x \otimes \ud x +|x|^{-2\alpha} \ud y \otimes \ud y\,.
\end{equation}
  The standard `\emph{Grushin plane}' corresponds to $\alpha=1$. Also for a Grushin-type plane one builds the Hilbert space $\cH_\alpha$, defined as in \eqref{eq:Hspacealpha}, and the Laplace-Beltrami differential operator $\Delta_{\mu_\alpha}:=\mathrm{div}_{\mu_\alpha}\circ\nabla$, explicitly given again by the analogue of 
   \eqref{eq:IntroLaplaceBeltrami}, and minimally realised as the analogue of \eqref{Halpha} on smooth functions compactly supported within each open half-plane. This yields the densely defined, non-negative, symmetric operator $H_\alpha$, and poses the problem of self-adjointness of $H_\alpha$, in order to analyse the generated heat or Schr\"{o}dinger flow.

   \begin{theorem}[\cite{Franceschi-Prandi-Rizzi-2017,Pozzoli_MSc2018,GMP-Grushin-2018,PozzoliGru-2020volume}]\label{thm:V-plane-conf-noconf} 
   %For $\alpha\geqslant 0$ let $M_\alpha\equiv(M,\mu_\alpha)$ and  be the Grushin plane introduced in \eqref{eq:Mgalpha} and let $H_\alpha$ be the operator on $\cH_\alpha=L^2(M,\ud\mu_\alpha)$ introduced in \eqref{eq:V-Halpha}.
 \begin{enumerate}
  \item[(i)] If $\alpha\in[-1,1)$, then $H_\alpha$ is not essentially self-adjoint in $\cH_\alpha$ and has infinite deficiency index.
  \item[(ii)] If $\alpha\in(-\infty,-1)\cup[1,+\infty)$, then $H_\alpha$ is essentially self-adjoint and therefore the Grushin-type plane $M_\alpha$ induces geometric quantum confinement.
 \end{enumerate}
\end{theorem}

  The above regime of essential self-adjointness was implicitly established in \cite{Franceschi-Prandi-Rizzi-2017} as an adaptation of the previous perturbative analysis \cite{Prandi-Rizzi-Seri-2016} devised for the compactified version of the manifold; the complete identification of essential self-adjointness and lack thereof was subsequently obtained in \cite{Pozzoli_MSc2018,GMP-Grushin-2018,PozzoliGru-2020volume} within a non-perturbative, novel scheme of \emph{constant-fibre direct integral decomposition} of the Hilbert space $\cH_\alpha=L^2(M,\ud\mu_\alpha)$ that generalises the direct integral decomposition \eqref{eq:Hxispace}-\eqref{eq:defF2} one performs in the compact case. This replaces uniformly fibred extensions on cylinder of the form \eqref{eq:unitequiv2} discussed above, namely,
  \[
   \bigoplus_{k\in\mathbb{Z}} A_{\alpha}(k)\,,
  \]
   $A_{\alpha}(k)$ acting self-adjointly on the fibre Hilbert space $\mathfrak{h}=L^2(\mathbb{R})$, with uniformly fibred direct integral extensions
  \[
   \int_{\mathbb{R}}^{\oplus} A_\alpha(\xi)\,\ud\xi\,,
  \]
  where the fibre operator $A_\alpha(\xi)$ on $\mathfrak{h}$ now depends on the continuous Fourier mode $\xi$, dual to the non-compact variable $y$.

  It is worth observing that the regime of self-adjointness for $\alpha$-Grushin cylinders and planes differ when $\alpha\in(-3,-1)$ (compare Theorems \ref{thm:cylinder} and \ref{thm:V-plane-conf-noconf}). This is due to the different nature of the direct sum and direct integral decompositions: indeed, when $\alpha\in(-3,-1)$, the only Fourier mode that is \emph{not} self-adjoint is the zero-th one, which brings a non-trivial contribution to the sum, but not to the integral. As a consequence, when $\alpha\in(-3,-1)$ the zero mode of a generic function $\psi\equiv\psi(x,y)$ hitting $\mathcal{Z}$ in the cylinder, namely the average on $\mathbb{S}^1$
  $$\psi_0(x)=\frac{1}{\sqrt{2\pi}}\int_{\mathbb{S}^1}\psi(x,y) \,\ud y\,,$$
  \emph{does} cross the singularity, whereas the zero mode of $\psi$ on the plane, namely
    $$\psi_0(x)=\frac{1}{\sqrt{2\pi}}\int_{\mathbb{R}}\psi(x,y) \,\ud y\,,$$
  \emph{does not} cross the singularity. The case $\alpha=-1$ is different as well between cylinder and plane: indeed, the non-self-adjoint Fourier modes are $\xi\in(-1,1)$ for the plane, and $k=0$ for the cylinder, thus yielding deficiency index of $H_\alpha$ equal to infinity for the plane, and equal to 2 for the cylinder.

  As a matter of fact, the lack of compactness makes the systematic identification of non-trivial self-adjoint extensions of $H_\alpha$ considerably harder and so far no explicit classification is available that mirrors Theorem \ref{thm:H_alpha_fibred_extensions} for the plane.
  %: only distinguished extensions (Friedrichs, Neumann, bridging) are known {\color{magenta} (references?)}.

  Beside the above concrete cylindrical and planar settings, the deep connection between geometry and self-adjointness is investigated for the problem of geometric confinement on more general almost-Riemannian structures. This includes `\emph{two-step two-dimensional almost-Riemannian structures}', characterised by an orthonormal frame for the metric in the vicinity of the singularity locus $\mathcal{Z}$ of the form \cite{Agrachev-Boscain-Sigalotti-2008}
  \begin{equation}
   X_1(x,y)\;=\;\dfrac{\partial}{\partial x}\,, \qquad X_2(x,y)\;=\;xe^{\phi(x,y)}\dfrac{\partial}{\partial y}
  \end{equation}
  (to be compared to \eqref{eq:frame} with $\alpha=1$). The essential self-adjointness of the corresponding minimally defined Laplace-Beltrami in the case of compactified $\mathcal{Z}$ was established in \cite{Boscain-Laurent-2013}.
%   , a work that also posed various further questions that are still open today
%     together with many still open questions.

   From a related perspective, the already observed circumstance that Grushin-type cylinders or planes are, classically, geodesically incomplete, but can induce, quantum-mechanically, geometric confinement (a condition that occurs more generally for regular almost-Riemannian manifold with compact singular set), poses an intriguing problem as far as semi-classical analysis is concerned. Indeed, reinstating Planck's constant in the Schr\"{o}dinger equation
   \begin{equation}\label{eq:Schreps}
    \ii\partial_t\psi+\varepsilon^2\Delta_{\mu_\alpha}\psi\;=\;0\,,\qquad\varepsilon>0
   \end{equation}
   (in the regime of $\alpha$ in which the minimally defined $\Delta_{\mu_\alpha}$ is unambiguously realised self-adjointly), semi-classics show, informally speaking, that as $\varepsilon\downarrow 0$ solutions get concentrated and evolves around geodesics. Therefore, the above-mentioned classical/quantum discrepancy makes the semi-classical analysis necessarily brake down in the limit.

   Such a discordance between classical and quantum picture can be at least partially resolved by appealing to different quantisation procedures on the considered Riemannian manifold, in practice considering corrections of the Laplace-Beltrami operator that have a suitable interpretation of free kinetic energy, much in the original spirit of \cite{DeWitt-1957}. Most of coordinate-invariant quantisation procedures (including path integral quantisation, covariant Weyl quantisation, geometric quantisation, and finite-dim\-en\-sional approximation to Wiener measures) modify $\Delta_{\mu_\alpha}$ with a term that depends on the scalar curvature $R_\alpha$ (which, in two dimensions, is twice the Gaussian curvature $K_\alpha$). This produces a replacement in \eqref{eq:Schreps}, in two dimensions, of $-\Delta_{\mu_\alpha}$ with the `\emph{curvature Laplacian}'
   \begin{equation}
    -\Delta_{\mu_\alpha}+cK_\alpha
   \end{equation}
   for suitable $c\geqslant 0$. In the recent work \cite{Boscain-Beschastnnyi-Pozzoli-2020} it was indeed shown, for generic two-step two-dimensional almost-Riemannian manifolds with compact singular set, that irrespective of $c\in(0,\frac{1}{2})$ the above correction washes essential self-adjointness out, yielding a quantum picture where the Schr\"{o}dinger particle does reach the singularity much as the classical particle does. (At the expenses of some further technicalities, the whole regime $c>0$ can be covered as well.) For concreteness, in the Grushin cylinder the effect of the curvature correction is evidently understood as a compensation between 
   $K=-\frac{2}{x^2}$ (see \eqref{eq:curvature} above) and the singular term $\frac{3}{4x^2}$ of the (unitary equivalent) Laplace-Beltrami operator. Still, the classical/quantum discrepancy discussed so far remains unexplained in more general settings.

   Concerning, instead, the heat flow, a satisfactory interpretation of the heat-confinement in the Grushin cylinder is known in terms of Brownian motions \cite{Boscain-Neel-2019} and random walks \cite{Agra-Bosca-Neel-Rizzi-2018}: roughly speaking, random particles are lost in the infinite area accumulated along $\mathcal{Z}$: the latter, in practice, acts as a barrier. Clearly, whereas curvature Laplacians are meaningful in the above context of inducing a non-confining (transmitting) Schr\"{o}dinger flow on two-step two-dimensional almost-Riemannian manifolds (including the Grushin cylinders), thus making quantum and classical picture more alike and well connected by semi-classics, this has no direct meaning instead in application to the heat flow on Riemannian or almost Riemannian manifolds. Indeed, as long as one regards the heat equation on manifold as a limit of a space-time discretised random walk, the stochastic process' generator \emph{is} the Laplace-Beltrami operator.
   
   %Introducing a correction (i.e., a rate of killing for the random particle) proportional to the curvature, the corresponding limiting random walk would unphysically fade away proportionally to the curvature. A 

  Generalisations of \cite{Boscain-Laurent-2013} have been established in \cite{Franceschi-Prandi-Rizzi-2017,Prandi-Rizzi-Seri-2016} from two-step two-dimensional almost-Riemannian structures to any dimensions, any step, and even to sub-Riemannian geometries, provided that certain geometrical assumptions on the singular set are taken. The main difficulty is the treatment of the `\emph{tangency}' (or `\emph{characteristic}') points: these are points belonging to the singularity of the metric structure where the vector distribution is tangent to the singularity. They are never present in Grushin cylinder or two-step almost-Riemannian structures, but may appear for example in three-step structures, such as, for instance,
\begin{equation}\label{ex:parabola}
X_1(x,y)\;=\;\dfrac{\partial}{\partial x}\,, \qquad X_2(x,y)\;=\;(y-x^2) \dfrac{\partial}{\partial y}\,,,\quad (x,y)\in \mathbb{R}^2 \,,
\end{equation}
 where the singularity is the parabola $y=x^2$ and the origin $(0,0)$ is a tangency point. Virtually nothing in known on the heat or the quantum confinement on such singular structures, including the simplest example \eqref{ex:parabola} (see \cite{Franceschi-Prandi-Rizzi-2015} for further remarks). First preliminary results in this respect were recently obtained in \cite{IB-2021}, where the interpretation of almost-Riemannian structures as special Lie manifolds permits to study some closure properties of singular perturbations of the Laplace-Beltrami operator even in the presence of tangency points. This opens new perspectives of treating several types of different singularities in sub-Riemannian geometry within the same unifying theory.

\section{A numerical glance at the bridging heat evolution}\label{sec:num}

 In this final Section we present and comment on qualitative features of the solution to the one-dimensional problem \eqref{eq:IVPgeneral}, \emph{obtained by numerical integration}, also in comparison with the initial value problem for the classical heat equation on $\mathbb{R}$.

 As already argued, it is the determination of the (integral kernel of) the heat propagator $\exp(-tA_\alpha^{\mathrm{B}})$, $t>0$, to be hard analytically, and this is due to the presence of boundary conditions for the solution at $x=0$ and any positive time.

 Numerics then represent a first, valuable way to access relevant aspects of the transmission of the heat flow between positive and negative half-line with bridging boundary conditions, and one may envisage that a systematic comparison will be launched numerically between analogous heat flows with different transmission protocols among those surveyed in Section \ref{sec:Grucyl}. Ours, here, is only an initial numerical glance at the bridging heat evolution to provide some insight and anticipate future investigations.

 Our numerical approach consists in approximating the solution $u=e^{-tA_\alpha^{\mathrm{B}}}\varphi$ to the problem \eqref{eq:IVPgeneral} by means of an approximated version of both the \emph{spatial} convolution integral between propagator's kernel and $\varphi$, and the \emph{complex} line integral that turns the resolvent of $A_\alpha^{\mathrm{B}}$ into its semi-group.

 More precisely, let us write
 \begin{equation}\label{eq:ukerphi}
  u(t,x)\;=\;(e^{-tA_\alpha^{\mathrm{B}}}\varphi)(x)\;=\;\int_{\mathbb{R}}\mathcal{K}_{\alpha}^{\mathrm{B}}(t;x,y)\varphi(x)\,\ud x
 \end{equation}
 where $\mathcal{K}_\alpha^{\mathrm{B}}(\cdot,\cdot)$ is the integral kernel of the bridging heat propagator. In turn, let us exploit the relation
 \begin{equation}\label{eq:invLapl}
  \begin{split}
   e^{-tA_\alpha^{\mathrm{B}}}\;&=\;\mathscr{L}^{-1}\big( (A_\alpha^{\mathrm{B}} - (\cdot) \mathbbm{1})^{-1}\big)(t) \\
   &=\;\frac{1}{2\pi \ii}\int_{\Gamma} e^{-z t}\big( (A_\alpha^{\mathrm{B}} - z \mathbbm{1})^{-1}\big)\,\ud z\,,\qquad t>0\,,
  \end{split}
 \end{equation}
 as an identity between bounded operators on $L^2(\mathbb{R})$ and with the integral understood in the Riemann sense in the strong operator topology, $\Gamma$ being a straight line in $\mathbb{C}$ orthogonal to the real axis in the open left half-plane, and $\mathscr{L}^{-1}$ denoting the inverse Laplace transform (the non-negativity of $A_\alpha^{\mathrm{B}}$ has led here to the choice $\mathfrak{Re}z<0$). \eqref{eq:invLapl} connects the resolvent of $A_\alpha^{\mathrm{B}}$ at the complex point $z$ with the semi-group at time $t>0$, and in terms of the integral kernels $(A_\alpha^{\mathrm{B}} - z \mathbbm{1})^{-1}(x,y)$ of the resolvent and $\mathcal{K}_{\alpha}^{\mathrm{B}}(t;x,y)$ of the propagator it reads
 \begin{equation}\label{eq:kerKalpha}
  \mathcal{K}_{\alpha}^{\mathrm{B}}(t;x,y)\;=\;\frac{1}{2\pi \ii}\int_{-1-\ii\cdot\infty}^{-1+\ii\cdot\infty} e^{-z t}\big( (A_\alpha^{\mathrm{B}} - z \mathbbm{1})^{-1}\big)(x,y)\,\ud z\,,\qquad t>0\,.
 \end{equation}
  The combinations of \eqref{eq:ukerphi} and \eqref{eq:kerKalpha} produces the solution $u$ and the two integrations contained therein may be computed numerically with standard packages.

  Of course, for \eqref{eq:ukerphi} and \eqref{eq:kerKalpha} to be implementable one needs to know the (integral kernel of) the resolvent $(A_\alpha^{\mathrm{B}} - z \mathbbm{1})^{-1}$. This is a not so hard knowledge to achieve from the underlying structure \eqref{eq:defineAbridge1d} of the operator $A_\alpha^{\mathrm{B}}$, once it is interpreted as a self-adjoint extension of the differential operator \eqref{eq:diffop} minimally defined on smooth functions compactly supported on $\mathbb{R}$ away from the origin. For this status of extension operator, one can appeal to the general Kre{\u\i}n-Vi\v{s}ik-Birman theory of self-adjoint extensions of lower semi-bounded and densely defined symmetric operators on Hilbert space \cite{GMO-KVB2017}, and obtain $(A_\alpha^{\mathrm{B}} - z \mathbbm{1})^{-1}$ fairly explicitly.

  The net result of this computation gives the following expression for the integral kernel %$\mathcal{R}_\alpha^{\mathrm{B}}(z)$ 
  of $(A_\alpha^{\mathrm{B}} - z \mathbbm{1})^{-1}$. With respect to the canonical decomposition
  \begin{equation}
   L^2(\mathbb{R},\ud x)\;\xrightarrow{\;\cong\;}\;L^2(\mathbb{R}^+,\ud x)\oplus L^2(\mathbb{R}^-,\ud x)\,,\qquad u\mapsto \begin{pmatrix} u^+ \\ u^- \end{pmatrix}
  \end{equation}
  (that is, $u^{\pm}(x):=u(x)$ for $x\gtrless 0$), consider the unitary transformation
  \begin{equation} 
  \begin{split}
   & U \,: \,L^2(\mathbb{R}^+,\ud x)\oplus L^2(\mathbb{R}^-,\ud x)\;\xrightarrow{\;\cong\;}\;L^2(\mathbb{R}^+,\ud x)\oplus L^2(\mathbb{R}^+,\ud x) \\
   & \qquad \qquad U\begin{pmatrix} u^+ \\ u^- \end{pmatrix}(x)\;=\;\begin{pmatrix} u^+(x) \\ u^-(-x) \end{pmatrix}\,,\qquad x>0\,,
  \end{split}
  \end{equation}
  and set $ \mathcal{R}_\alpha^{\mathrm{B}}(z):=U(A_\alpha^{\mathrm{B}} - z \mathbbm{1})^{-1} U^{-1}$. Then
  \begin{equation}\label{eq:resolventPsandwitched}
  (A_\alpha^{\mathrm{B}} - z \mathbbm{1})^{-1}\;=\;U^{-1}\,\mathcal{R}_\alpha^{\mathrm{B}}(z)\, U
  \end{equation}
   and the integral kernel of $\mathcal{R}_\alpha^{\mathrm{B}}(z)$ is given by
  \begin{equation}\label{eq:Rkernel}
  \begin{split}
   \mathcal{R}_\alpha^{\mathrm{B}}(z)(x,y)\;&=\;\mathcal{G}_{\alpha,z}(x,y)\begin{pmatrix} 1 & 0 \\ 0 & 1 \end{pmatrix}-\frac{\ii\pi}{8}\cos\Big(\frac{\pi\alpha}{2}\Big)\, e^{\ii\frac{\pi\alpha}{2}}\begin{pmatrix} 1 & 1 \\ 1 & 1 \end{pmatrix}P_{\alpha,z}(x)P_{\alpha,z}(y)\,, \\
   & \qquad x>0\,,\; y>0\,,
  \end{split}
  \end{equation}
 where, in terms of the Bessel functions of first and second kind $J_\nu$ and $Y_\nu$,
 \begin{equation}
  \begin{split}
   P_{\alpha,z}(x)\;&=\;\sqrt{x}\,J_{\frac{1+\alpha}{2}}(x\sqrt{z})+\ii\sqrt{x}\,Y_{\frac{1+\alpha}{2}}(x\sqrt{z}) \\
   Q_{\alpha,z}(x)\;&=\;2\sqrt{x}\,J_{\frac{1+\alpha}{2}}(x\sqrt{z})
  \end{split}\qquad (\mathfrak{Im}\sqrt{z}>0)\,,
 \end{equation}
 and
 \begin{equation}
  \mathcal{G}_{\alpha,z}(x,y)\;=\;
  -\frac{\ii\pi}{4}\begin{cases}
   \; P_{\alpha,z}(x)\,Q_{\alpha,z}(y)\,, & \textrm{ if }\;\; 0<y<x\,, \\
   \; Q_{\alpha,z}(x)\,P_{\alpha,z}(y)\,, & \textrm{ if }\;\; 0<x<y\,.
  \end{cases}
 \end{equation}

       \begin{figure}[t!]
     \begin{center}
      \includegraphics[width=9.8cm]{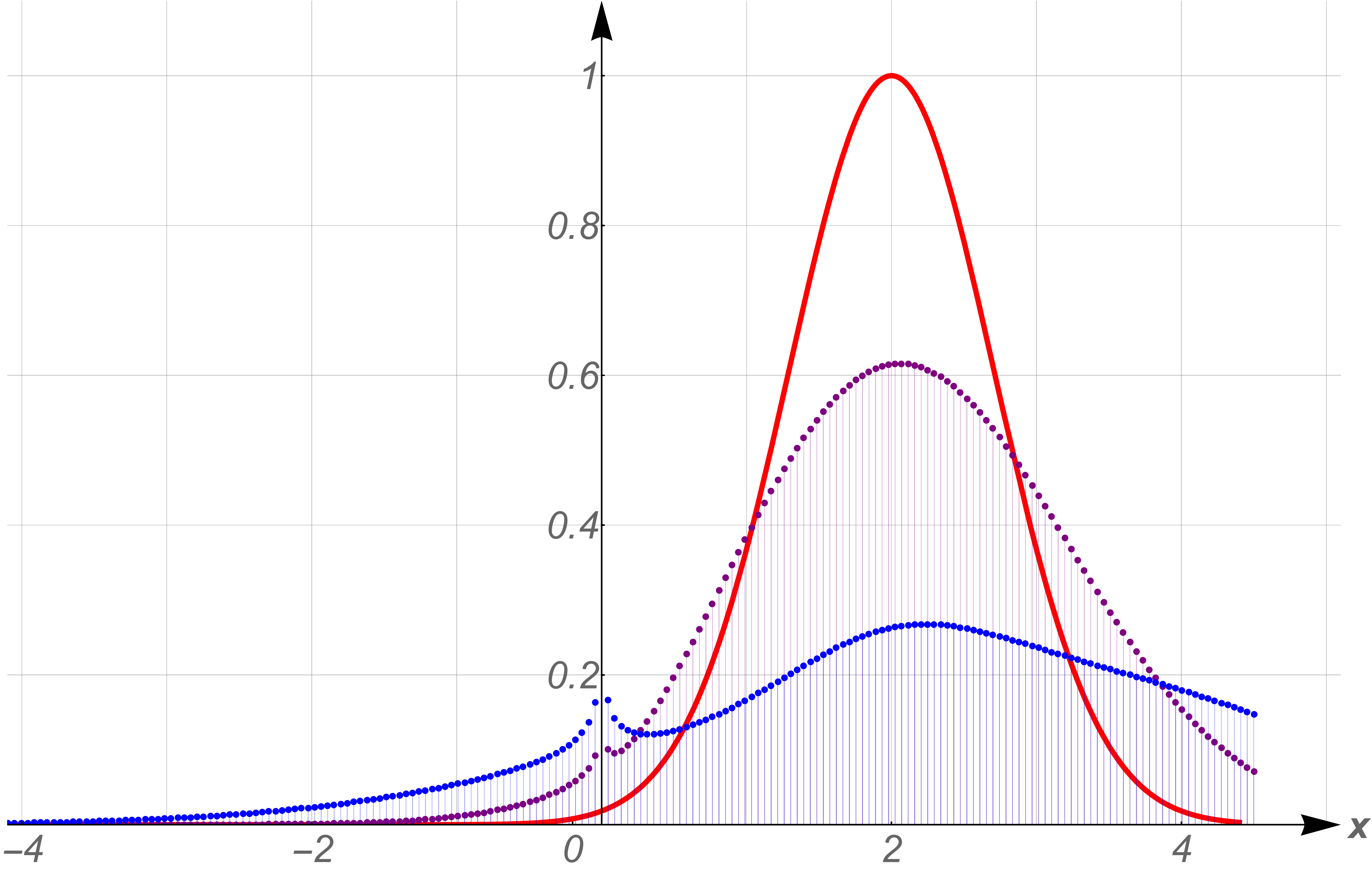}\caption{Solution $u(t,x)$ to the heat-bridging initial value problem \eqref{eq:IVPgeneral} with Gaussian initial datum $\varphi(x)=e^{-(x-2)^2}$ (red curve). Plot of $|u(t,\cdot)|$ at $t=0.5$ (magenta dotted line) and $t=2$ (blue dotted line).\label{fig:multigau}}
     \end{center}
\end{figure}

      \begin{figure}[t!]
     \begin{center}
      \includegraphics[width=9.8cm]{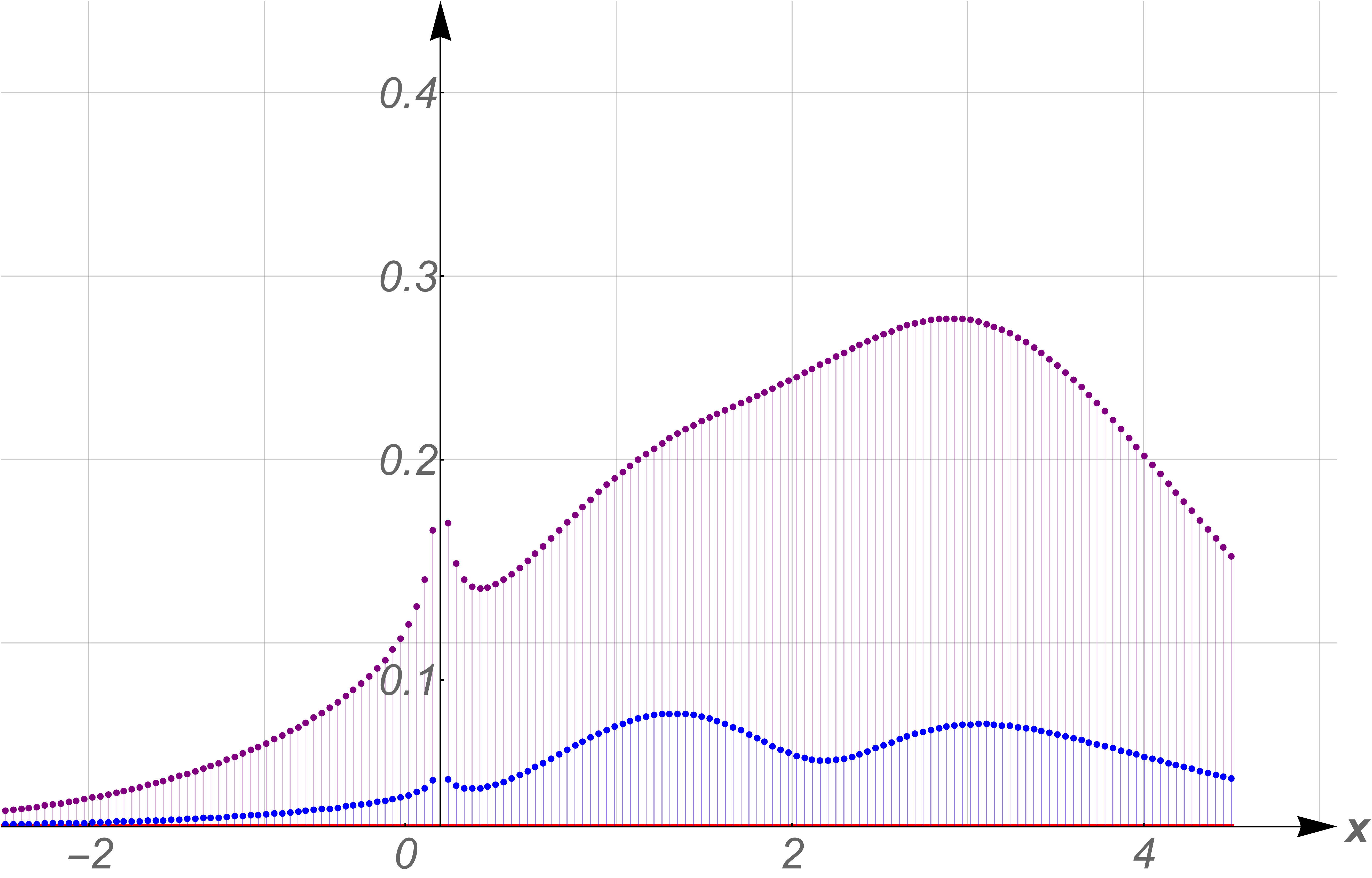}\caption{Comparison at time $t=1.5$ between the solution to the heat-bridging initial value problem \eqref{eq:IVPgeneral} with zero-momentum Gaussian initial datum $\varphi(x)=e^{-(x-2)^2}$ (magenta dotted curve) and with non-zero momentum Gaussian initial datum $\varphi_2(x)=e^{-3\ii x}e^{-(x-2)^2}$ towards left. Both plots are of $|u(t,\cdot)|$. The evolution of the Gaussian initially shot towards left displays in-going + outgoing oscillation.\label{fig:gauoscill}}
     \end{center}
\end{figure}

      \begin{figure}[t!]
     \begin{center}
      \includegraphics[width=9.8cm]{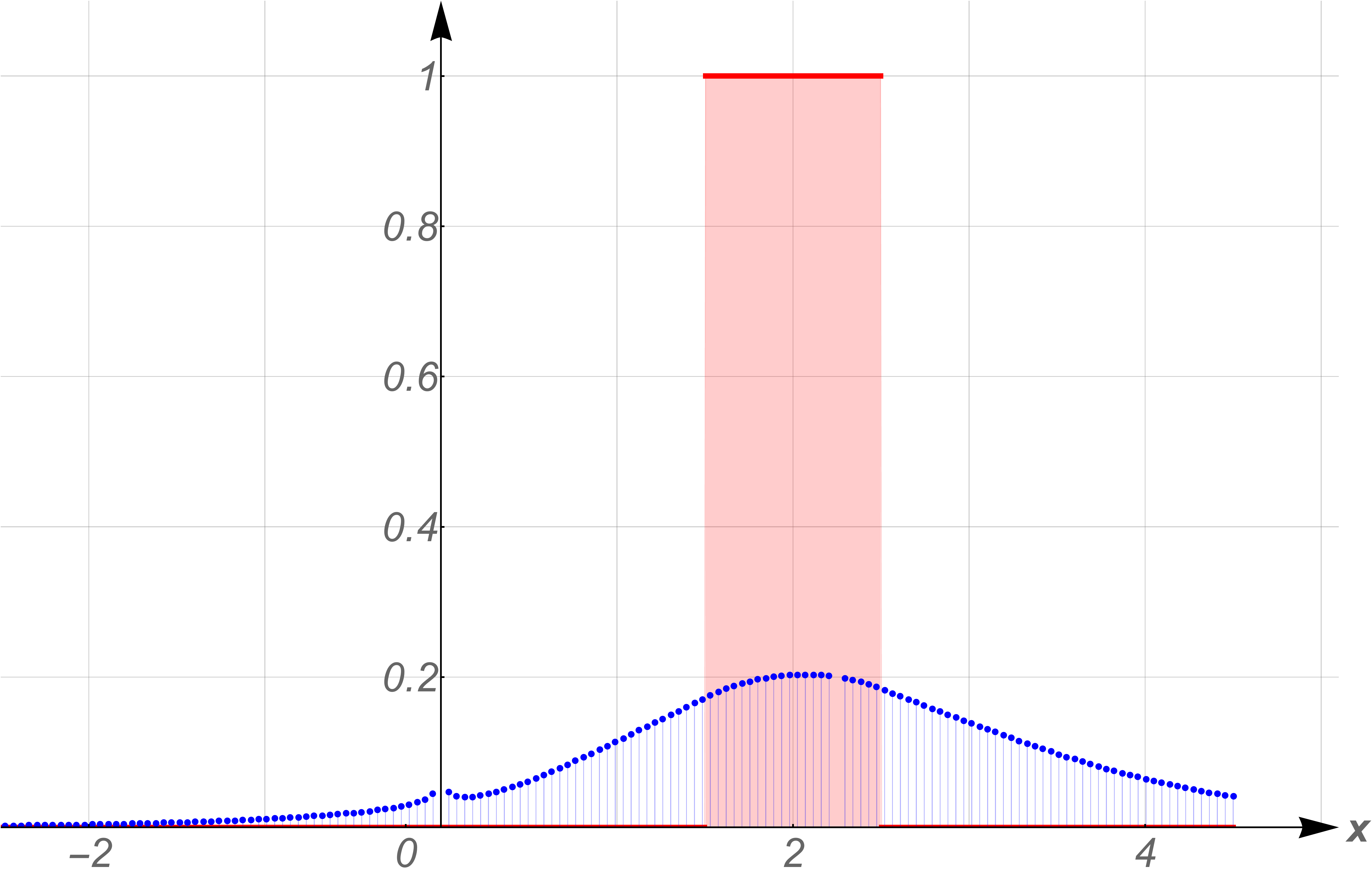}\caption{Solution $u(t,x)$ to the heat-bridging initial value problem \eqref{eq:IVPgeneral} with initial datum $\varphi(x)$ given by the characteristic function of the interval $[\frac{1}{2},\frac{3}{2}]$ (red curve). Plot of $|u(t,\cdot)|$ at $t=0.5$ (blue dotted line)\label{fig:stepev}}
     \end{center}
\end{figure}

        \begin{figure}[h!]
     \begin{center}
      \includegraphics[width=9.8cm]{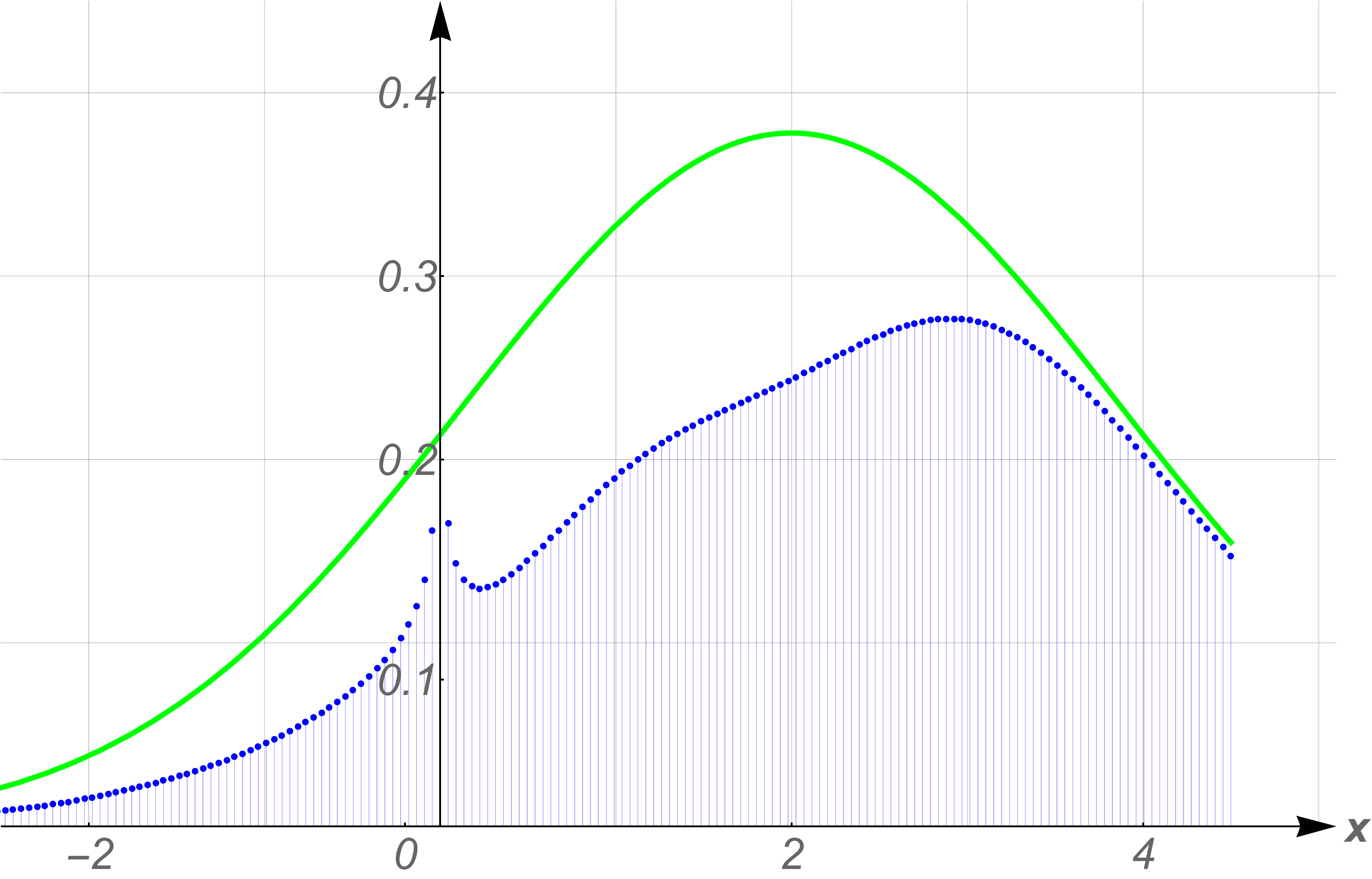}\caption{Comparison at time $t=1.5$ between the solution to the heat-bridging initial value problem \eqref{eq:IVPgeneral} with Gaussian initial datum $\varphi(x)=e^{-(x-2)^2}$ (blue dotted curve) and the solution to the ordinary heat equation on $\mathbb{R}$ (green curve).\label{fig:freecal}}
     \end{center}
\end{figure}

  When formulas \eqref{eq:ukerphi}, \eqref{eq:kerKalpha}, \eqref{eq:resolventPsandwitched}, and \eqref{eq:Rkernel} are implemented numerically we obtain a scenario exemplified in Figures \ref{fig:multigau}, \ref{fig:gauoscill}, \ref{fig:stepev}, \ref{fig:freecal} below.

  For concreteness, the bridging-heat evolution is considered, namely the solution $u\equiv u(t,x)$ to \eqref{eq:IVPgeneral}, of an initial datum $\varphi$ essentially supported on the right half-line. An initial Gaussian is seen to evolve at later times with the typical heat-flow flattening of the solution, with the immediate formation of the characteristic bridging behaviour at $x=0$ (Fig.~\ref{fig:multigau}).

  Notably, if $\varphi$ is additionally shot with an initial non-zero momentum towards the singularity, its evolution displays an oscillation given by the superposition of an in-going wave and a component that bounces backwards (Fig.~\ref{fig:gauoscill}), as compared with the evolution at the same time of the same Gaussian with no initial momentum.

  It is also pretty transparent that the bridging-heat flow has a regularising effect at every $t>0$, as observed with the evolution of an initial step function. (Fig.~\ref{fig:stepev}).

  We have also further evidence of a qualitatively similar behaviour of the free heat flow and the bridging-heat flow, but of course for the characteristic boundary condition of bridging type at the origin (Fig.~\ref{fig:freecal}).

  Whereas, as said, this provides only a first glance at the qualitative properties of the bridging-heat evolution on two connected half-lines, the evidences collected here are encouraging and further corroborate the quest for the analytic identification of counterpart $L^p$-$L^q$ estimates, smoothing estimates, and space-time (Strichartz) estimates for the bridging-heat flow, as compared to \eqref{eq:heatdispest} for the classical heat flow.

\section*{Acknowledgement}
This work is partially supported by the QUACO project (grant ANR-17-CE-40-0007-01), the EIPHI Graduate School (ANR-17-EURE-0002), the CONSTAT project (funded by the Conseil R\'egional de Bourgogne Franche Comt\'e and the European Union through the PO FEDER Bourgogne 2014/2020 programmes), the MIUR-PRIN 2017 project MaQuMA cod.~2017ASFLJR, the INdAM -- Italian National Institute for Higher Mathematics, and the Alexander von Humboldt Foundation.

 \def\cprime{$'$}


\begin{thebibliography}{10}

\bibitem{Agra-Bosca-Neel-Rizzi-2018}
{\sc A.~Agrachev, U.~Boscain, R.~Neel, and L.~Rizzi}, {\em {Intrinsic random
  walks in {R}iemannian and sub-{R}iemannian geometry {\it via} volume
  sampling}}, ESAIM Control Optim. Calc. Var., 24 (2018), pp.~1075--1105.

\bibitem{Agrachev-Boscain-Sigalotti-2008}
{\sc A.~Agrachev, U.~Boscain, and M.~Sigalotti}, {\em {A {G}auss-{B}onnet-like
  formula on two-dimensional almost-{R}iemannian manifolds}}, Discrete Contin.
  Dyn. Syst., 20 (2008), pp.~801--822.

\bibitem{IB-2021}
{\sc I.~Beschastnyi}, {\em {Closure of the Laplace-Beltrami operator on 2D
  almost-Riemannian manifolds and semi-Fredholm properties of differential
  operators on Lie manifolds}}, arXiv:2104.07745 (2021).

\bibitem{Boscain-Beschastnnyi-Pozzoli-2020}
{\sc I.~Beschastnyi, U.~Boscain, and E.~Pozzoli}, {\em {Quantum Confinement for
  the Curvature Laplacian $-\Delta + cK$ on 2D-Almost-Riemannian Manifolds}},
  Potential Analysis,  (2021).

\bibitem{Boscain-Laurent-2013}
{\sc U.~Boscain and C.~Laurent}, {\em {The {L}aplace-{B}eltrami operator in
  almost-{R}iemannian geometry}}, Ann. Inst. Fourier (Grenoble), 63 (2013),
  pp.~1739--1770.

\bibitem{Boscain-Neel-2019}
{\sc U.~Boscain and R.~W. Neel}, {\em {Extensions of {B}rownian motion to a
  family of {G}rushin-type singularities}}, Electron. Commun. Probab., 25
  (2020), pp.~Paper No. 29, 12.

\bibitem{Boscain-Prandi-JDE-2016}
{\sc U.~Boscain and D.~Prandi}, {\em {Self-adjoint extensions and stochastic
  completeness of the {L}aplace-{B}eltrami operator on conic and anticonic
  surfaces}}, J. Differential Equations, 260 (2016), pp.~3234--3269.

\bibitem{Boscain-Prandi-Seri-2014-CPDE2016}
{\sc U.~Boscain, D.~Prandi, and M.~Seri}, {\em {Spectral analysis and the
  {A}haronov-{B}ohm effect on certain almost-{R}iemannian manifolds}}, Comm.
  Partial Differential Equations, 41 (2016), pp.~32--50.

\bibitem{Burq-Planchon-Stalker-2003}
{\sc N.~Burq, F.~Planchon, J.~G. Stalker, and A.~S. Tahvildar-Zadeh}, {\em
  {Strichartz estimates for the wave and {S}chr{\"o}dinger equations with the
  inverse-square potential}}, J. Funct. Anal., 203 (2003), pp.~519--549.

\bibitem{Calin-Chang-SubRiemannianGeometry}
{\sc O.~Calin and D.-C. Chang}, {\em {Sub-{R}iemannian geometry}}, vol.~126 of
  {Encyclopedia of Mathematics and its Applications}, Cambridge University
  Press, Cambridge, 2009.
\newblock General theory and examples.

\bibitem{Derezinski-Georgescu-2021}
{\sc J.~Derezi{\'n}ski and V.~Georgescu}, {\em {On the domains of {B}essel
  operators}}, Ann. Henri Poincar{\'e}, 22 (2021), pp.~3291--3309.

\bibitem{DeWitt-1957}
{\sc B.~S. DeWitt}, {\em {Dynamical Theory in Curved Spaces. I. A Review of the
  Classical and Quantum Action Principles}}, Rev. Mod. Phys., 29 (1957),
  pp.~377--397.

\bibitem{Franceschi-Prandi-Rizzi-2015}
{\sc V.~Franceschi, D.~Prandi, and L.~Rizzi}, {\em {Recent results on the
  essential self-adjointness of {sub-Laplacians,} with some remarks on the
  presence of characteristic points}}, S{\'e}minaire de th{\'e}orie spectrale
  et g{\'e}om{\'e}trie, 33 (2015-2016), pp.~1--15.

\bibitem{Franceschi-Prandi-Rizzi-2017}
\leavevmode\vrule height 2pt depth -1.6pt width 23pt, {\em {On the essential
  self-adjointness of singular sub-{L}aplacians}}, Potential Anal., 53 (2020),
  pp.~89--112.

\bibitem{GM-Grushin3-2020}
{\sc M.~Gallone and A.~Michelangeli}, {\em {Quantum particle across Grushin
  singularity}}, Journal of Physics A: Mathematical and Theoretical,  (2021).

\bibitem{GMO-KVB2017}
{\sc M.~Gallone, A.~Michelangeli, and A.~Ottolini}, {\em
  {Kre{\u\i}n-Vi\v{s}ik-Birman self-adjoint extension theory revisited}}, in
  {Mathematical Challenges of Zero Range Physics}, A.~Michelangeli, ed.,
  {INdAM-Springer series, Vol.~42}, Springer International Publishing, 2020,
  pp.~239--304.

\bibitem{GMP-Grushin-2018}
{\sc M.~Gallone, A.~Michelangeli, and E.~Pozzoli}, {\em {On geometric quantum
  confinement in {G}rushin-type manifolds}}, Z. Angew. Math. Phys., 70 (2019),
  pp.~Art. 158, 17.

\bibitem{GMP-Grushin2-2020}
\leavevmode\vrule height 2pt depth -1.6pt width 23pt, {\em {Geometric
  confinement and dynamical transmission of a quantum particle in Grushin
  cylinder}}, arXiv:2003.07128 (2020).

\bibitem{M-2015-nonStrichartzHartree}
{\sc A.~Michelangeli}, {\em {Global well-posedness of the magnetic {H}artree
  equation with non-{S}trichartz external fields}}, Nonlinearity, 28 (2015),
  pp.~2743--2765.

\bibitem{PozzoliGru-2020volume}
{\sc E.~Pozzoli}, {\em {Quantum Confinement in $\alpha$-Grushin planes}}, in
  {Mathematical Challenges of Zero-Range Physics}, A.~Michelangeli, ed.,
  {Springer INdAM Series}, Springer International Publishing, 2021,
  pp.~229--237.

\bibitem{Pozzoli_MSc2018}
\leavevmode\vrule height 2pt depth -1.6pt width 23pt, {\em {Models of quantum
  confinement and perturbative methods for point interactions}}, Master Thesis
  (2018).

\bibitem{Prandi-Rizzi-Seri-2016}
{\sc D.~Prandi, L.~Rizzi, and M.~Seri}, {\em {Quantum confinement on
  non-complete {R}iemannian manifolds}}, J. Spectr. Theory, 8 (2018),
  pp.~1221--1280.

\bibitem{rs2}
{\sc M.~Reed and B.~Simon}, {\em {Methods of modern mathematical physics. {II}.
  {F}ourier analysis, self-adjointness}}, Academic Press [Harcourt Brace
  Jovanovich, Publishers], New York-London, 1975.

\bibitem{schmu_unbdd_sa}
{\sc K.~Schm{\"u}dgen}, {\em {Unbounded self-adjoint operators on {H}ilbert
  space}}, vol.~265 of {Graduate Texts in Mathematics}, Springer, Dordrecht,
  2012.

\bibitem{Wang-Zhaohui-Chengchun_HarmonicAnalysisI_2011}
{\sc B.~Wang, Z.~Huo, C.~Hao, and Z.~Guo}, {\em {Harmonic analysis method for
  nonlinear evolution equations. {I}}}, World Scientific Publishing Co. Pte.
  Ltd., Hackensack, NJ, 2011.

\end{thebibliography}
\end{document}